\documentclass{cocv}
\usepackage{lineno,hyperref}
\modulolinenumbers[5]
\usepackage{mathrsfs}
\usepackage{amsfonts}
\usepackage{graphics}
\usepackage{graphicx}
\usepackage{mathtools,amsfonts,amssymb,amsthm,mathrsfs,setspace}
\usepackage{bm,verbatim}
\usepackage{subfig}
\usepackage{lipsum}
\usepackage[dvips]{epsfig}
\usepackage{verbatim}
\usepackage{mathrsfs}
\usepackage{hyperref}
\usepackage{indentfirst}
\usepackage{stmaryrd}
\usepackage{subeqnarray}
\usepackage{xcolor,pgfplots,float,layouts,booktabs,wrapfig}
\usepackage[capitalise]{cleveref}
\usepackage[absolute,overlay]{textpos}

\usepackage{xcolor,pgfplots,float,layouts,booktabs,wrapfig}
\usetikzlibrary{arrows}

\newtheorem{theorem}{Theorem}[section]
\newtheorem{lemma}[theorem]{Lemma}

\newtheorem{corollary}[theorem]{Corollary}
\newtheorem{remark}{Remark}[section]

\newcommand{\mc}{\mathcal}

\renewcommand{\thepage}{}

\def\e1{{\varepsilon_{11}}}
\def\b1{{\beta_{11}}}
\def\bp3{{\beta_{33}}}
\def\ep3{{\varepsilon_{33}}}

\def\Ltwo{{\mathbb L}^2 }

\newcommand{\R}{\mathbb{R}}

\newcommand{\N}{\mathbb{N}}

\usepackage{setspace}
\allowdisplaybreaks
\newcommand{\mathsym}[1]{{}}
\newcommand{\unicode}[1]{{}}

\begin{document}
\title{A novel sensor design for a cantilevered Mead-Marcus-type sandwich beam model by the order-reduction technique}\thanks{A. \"{O}. \"{O}zer  gratefully acknowledges the financial support of the National Science Foundation of USA under Cooperative Agreement No: 1849213.}
\author{ A.O. Ozer
}\address{Department of Mathematics, Western Kentucky University, Bowling Green, KY 42101, USA}
\author{A.K. Aydin}\address{Department of Mathematics, Western Kentucky University, Bowling Green, KY 42101, USA}

%
%
\begin{abstract}
{\color{black}A novel space-discretized Finite Differences-based model reduction introduced in \cite{Guo3} is extended to the partial differential equations (PDE) model of a multi- layer Mead-Marcus-type sandwich beam with clamped-free boundary conditions.}
The PDE model describes transverse vibrations for a sandwich beam whose alternating outer elastic layers constrain viscoelastic core layers, which allow transverse shear. The major goal of this project is to design a single tip velocity sensor to control the overall dynamics on the beam. Since the spectrum of the PDE can not be constructed analytically, the so-called multipliers approach is adopted to prove that the PDE model is exactly observable with sub-optimal observation time.  Next, the PDE model is reduced by the ``order-reduced'' Finite-Differences technique. This method does not require any type of filtering though the exact observability as $h\to 0$ is achieved by a constraint on the material constants. The main challenge here is the strong coupling of the shear dynamics of the middle layer with overall bending dynamics. This complicates the absorption of coupling terms in the discrete energy estimates. This is sharply different from a single-layer (Euler-Bernoulli) beam. \end{abstract}
%
%
\maketitle


\section{Introduction}\label{sec-intr}
  Multi-layer  ``sandwich'' beams  have become popular in industrial applications in aeronautic, civil, defense, biomedical, and space structures, e.g. see \cite{biomed,Wu}. A multi-layer sandwich beam is a layered structure consisting of perfectly bonded alternating elastic (stiff) layers  sandwiching  compliant viscoelastic layers \cite{Baz,Hansen3}. The viscoelastic layers are meant to allow the transverse shear so that the vibrations on the beam are passively suppressed. The most commonly-used sandwich beam theories in the literature avoid treating each layer separately by averaging in various ways the stresses and elastic moduli through the thickness. One of the widely-accepted discrete-layer theories treating each layer separately is  by Mead-Marcus where piezoelectric/elastic  layers adopt the Euler-Bernoulli small displacement assumptions and viscoelastic layers adopt the Mindlin-Timoshenko small displacement assumptions. Moreover, it is assumed that the longitudinal and rotational inertia terms of the outer layers are discarded, see  \cite{mead_forced_1969,Hansen3} for further discussions of models.

	\begin{figure}[htb!]
		\centering
		\includegraphics[width=14.5cm]{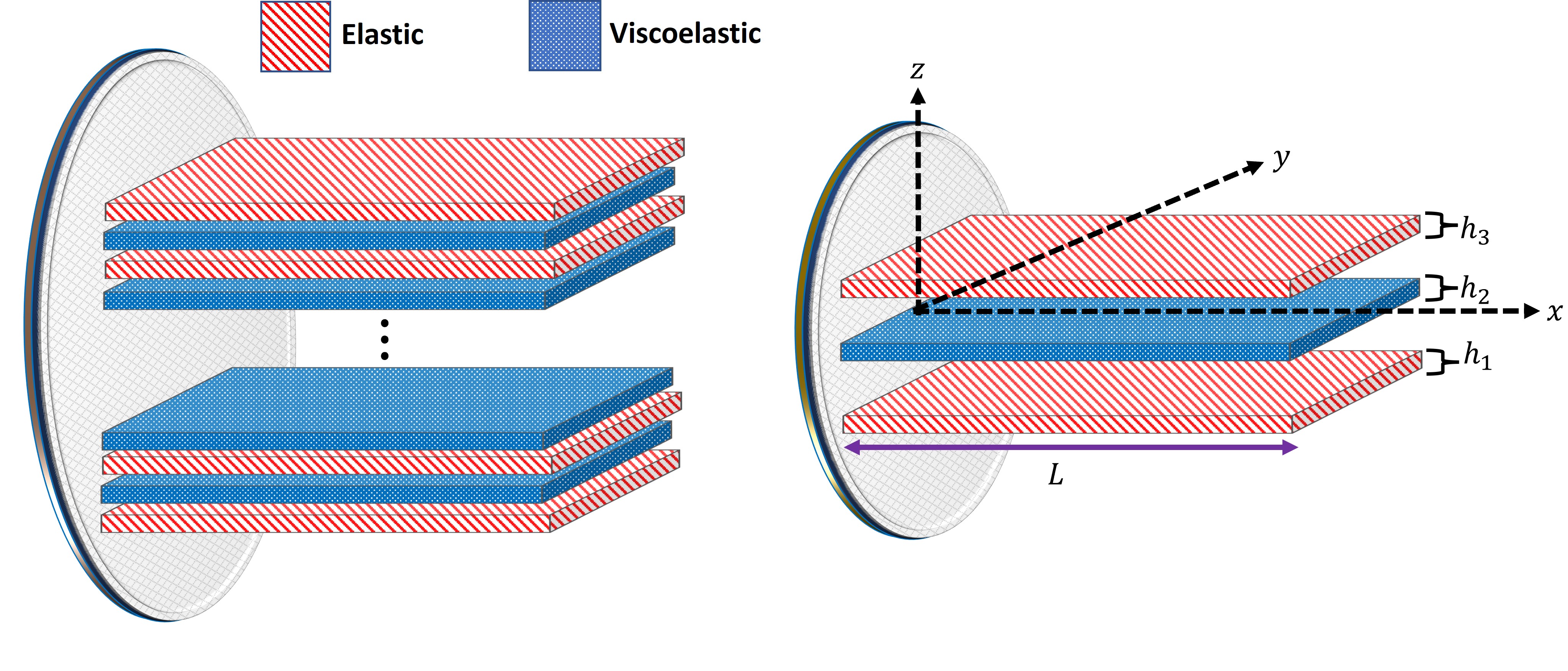}
		\caption{A multi-layer beam consisting of alternating elastic (i.e. aluminum, titanium, PZT) and viscoelastic layers (i.e. honeycomb) where the viscoelastic layers are constrained by the outer elastic layers (left). A three-layer beam of length $L$, in particular, consisting of two outer  layers (with thicknesses $h_1$ and $h_3$) constraining a viscoelastic layer (with thickness $h_2$) (right). }
		\label{f1}
	\end{figure}

 Letting $z(x,t)$  denote the uniform transverse (bending) displacement  of the centerlines of each layer and $\vec{\bm \phi}(x,t)$ denote an $m\times 1$ column vector of shear angles $v_i$ for the $i^{th}$ viscoelastic layer, see Fig. \ref{f1}, the equations of motion for an Mead-Marcus beam model  consisting with $(m+1)$ stiff layers  and $m$ viscoelastic layers are
	\begin{equation}\label{raymulti}  \left\{
		\begin{array}{ll}
			\ddot z + z'''' - {\bf B}^T  \vec{\bm \phi}' =0,&\\
			-{\bf C}  \vec{\bm \phi}'' + {\bf P}\vec{\bm \phi} = -Bz''',~~  (x,t)\in (0,L) \times \mathbb{R}^+&\\
			 z(0,t)=z'(0,t)=\vec{\bm \phi}(0,t)=z''(L,t)=z'''(L,t)-B \vec{\bm \phi}(L,t)=\vec{\bm \phi}'(L,t)=0,\qquad t \in \mathbb{R}^+ &\\
			(z,\dot z)(x,0)=(z_0,z_1)(x),~ x \in (0,L)
		\end{array}
		\right.
	\end{equation}
where ${\bf B}$ is a $m\times1$ column vector with positive entries, ${\bf C}$
and ${\bf P}$ are $ m\times m$ invertible, symmetric, positive definite matrices, {\color{black}  and they  are defined in terms of shear moduli, Young's moduli, Poisson's ratios, and thicknesses $h_i$ of the layers \cite{Hansen3,O-H-CDC}.  The three-layer version \cite{Fabiano-Hansen,ozer_modeling_2016,IFAC,ozer_dynamic_2018} of \eqref{raymulti} can be obtained by considering  $B,C,P$  positive constants (defined in terms of material constants) and $\vec {\bm \phi}=\phi$ (scalar-valued) describing the shear of the only middle layer as the following
	\begin{eqnarray}  \label{eq1a}\left\{
		\begin{array}{ll}
			\ddot{z}+z''''-B\phi'=0, & \\
			-C\phi''+P\phi=-Bz''',&(x,t)\in (0,L)\times\R^+
		\end{array}\right.
	\end{eqnarray}
	with clamped-free boundary conditions
	\begin{eqnarray} \label{eq1c}
		\begin{array}{ll}
			z(0,t)=z'(0,t)=\phi(0,t)=0,
			z''(L,t)=z'''(L,t)-B\phi(L,t)=\phi'(L,t)=0,& t\in  \R^+
		\end{array}
	\end{eqnarray}
	and initial conditions
	\begin{eqnarray}  \label{eq1d}
		\begin{array}{ll}
			z(x,0)=z^0(x),\quad \dot z(x,0)=z^1(x),& x\in[0,L].
		\end{array}
	\end{eqnarray}
	
	There are results in the literature for the controllability of the PDE models  \eqref{raymulti} and \eqref{eq1a} with fully-clamped boundary conditions and  with a single boundary controller \cite{O-H-CDC,rajeev_moment_2003}. The exponential stability  and the analyticity of   \eqref{raymulti} with hinged boundary condition are studied with viscous damping proportional to the shear in the viscoelastic layers  \cite{Allen-Hansen}. Recently, \eqref{eq1a} with  PZT outer layer  is shown to be exponentially stable with a single feedback controller \cite{IFAC} placed at the moment boundary condition.

	The literature for approximating PDE-based models for waves, beams, or plates suggests that the available computational techniques are not mathematically rigorous or  reliable if one implements them ``blindly'', i.e. see \cite{zuazua_controllability_2006}. This issue is even more pronounced for multi-layer beams where there are strong coupling of vibrations on the discrete layers. For example, a  highly popular  empirical technique, called the Finite Element  approximations \cite{A-A,B-V},  considers only the first several low-frequency vibrational modes due to its computational feasibility. However, if  only a few vibration modes are considered to design a controller/sensor and the high-frequency residual modes are completely ignored, the   so-called {spill-over} effect may occur whereby these unobserved (unsensed) modes could prevent observability/stabilization of the whole system \cite{Ballas,Hebrard,Y-W}. Another highly popular method is called the Finite Difference approximations \cite{B-M-Z,Micu,E-Z,Roventa}.  However, a blind  application of this approximation generates {spurious high-frequency} unobservable vibrational modes in the approximated dynamics \cite{I-Z}. These modes cannot be distinguished by the sensor design. Therefore, the finite dimensional counterpart do not satisfy the exact observability property uniformly as the discretization parameter approaches zero, whereas the infinite-dimensional system (PDE) is still exactly observable. The main reason for the lack of exact observability is the accumulation of high-frequency eigenvalues to a specific value as the discretization parameter approaches zero. These high-frequency eigenvalues can be eliminated by the so-called \textit{``direct Fourier filtering''} or/and the \textit{``indirect filtering''} \cite{T-Z} methods. In fact, this is successfully applied to the single Euler-Bernoulli beam  \cite{L-Z}, Rayleigh beam \cite{ozer_uniform_2019}, and hinged multi-layer beam  \cite{Aydin-O1,Aydin-LCSS}. The direct Fourier filtering method encompasses both low and high-frequency eigenvalues, while ignoring spurious eigenvalues.
	
 Recently,  a novel  semi-discretized numerical scheme  based on the  order-reduction technique
is reported for the sensor design of a single cantilevered beam \cite{Ren,Guo1,Guo2}. The proposed model reduction is constructed on equidistant grid points and average operators.  The uniform observability of the discretized model can be proved by the discrete multipliers method, very similarly to the continuous counterpart. This is particularly compatible with cantilevered beams since the Finite Difference-based model reductions, in general, do not allow the explicit description of their spectrums.

To the best of our knowledge, there are no results reported in the literature for the robust sensor designs for reduced models corresponding to \eqref{eq1a}-\eqref{eq1d}.
The major contribution of this work is that the exact observability of both the  PDE model  \eqref{eq1a} and its model reduction \eqref{clampeddiscrete} are rigorously established by the multipliers approach.  Our approach  not only extends the results in \cite{O-H-CDC,rajeev_moment_2003,IFAC} but also provides  better insights to understand the overall observability/sensor design of their reduced counterparts.

	\section{Well-posedness and Exact Observability Results}
	Consider the PDE model \eqref{eq1a}-\eqref{eq1d}. The second equation in \eqref{eq1a} is elliptic, and therefore,  $\phi$ can be solved and it can be written in terms of bending $z$. By defining the identity operator $D_x^0=I,$ and the differential operator $D_x^2:=\frac{d^2}{d x^2}$ on the domain ${\rm Dom}(D_x^2)=H^2_L(0,L):=\{\psi\in H^2(0,L): \psi_x(0)=\psi_x(L)=0\}.$ the operator $(-CD_x^2+PI),$ defined on the ${\rm Dom}(D_x^2)$, is positive definite and self-adjoint, and therefore, invertible. Therefore,
	\begin{eqnarray}
		&\phi=-B\left(-CD_x^2+PI\right)^{-1}D_x^2z'.\label{phi-def}
	\end{eqnarray}

\begin{lemma} \label{lemma1}
		Define $J:=-\frac{1}{C}+\frac{P}{C}(-CD_x^2+P I)^{-1}$. Then, $J$ is continuous, non-positive, and self-adjoint  on $L^2(0,L)$. Moreover, for all $u\in{\rm Dom}(D_x^2),$ $
			J=(-CD_x^2+P I)^{-1}D_x^2.$
	\end{lemma}
	\begin{proof} It is trivial to show that $J$ is continuous on $L^2(0,L)$. Now, let $w=(-CD_x^2+P I)^{-1}u$ so that $-CD_x^2w+P w=u$. Then, $J$ is non-positive since
		\begin{align*}
			<Ju,u>&=\left<-\frac{1}{C}u+\frac{P}{C}(-CD_x^2+P I)^{-1}u,u\right>
			=-C\left<D_x^2w,D_x^2w\right>-P\left<D_x w, D_x w\right>
			\le 0.
		\end{align*}
		Next, by letting $u,v\in L^2(0,L),$
		\begin{align*}
			\langle Ju,v\rangle &=\left\langle 	-\frac{1}{C}u+\frac{P}{C}(-CD_x^2+P I)^{-1}u,v\right\rangle= \left\langle 	u,-\frac{1}{C}v+\frac{P}{C}(-CD_x^2+P I)^{-1}v\right\rangle=\langle u,Jv\rangle
		\end{align*}
		which implies the self-adjoint-ness of $J.$
		Finally, let $Ju=-\frac{1}{C}u+\frac{P}{C}(-CD_x^2+P I)^{-1}u=v.$ Then, $u=\frac{-C^2D_x^2}{P}v+Cv-\frac{CD_x^2}{P}u+u,$ and therefore, $v=(-CD_x^2+P I)^{-1}D_x^2u.$
		Hence, $J=(-CD_x^2+P I)^{-1}D_x^2$.
	\end{proof}
	By the definition of $J$  in \cref{lemma1}, $\phi=-BJz',$ and the system \eqref{eq1a}-\eqref{eq1d} can be reformulated as
	\begin{eqnarray}  \label{compactclamped}\left\{
		\begin{array}{ll}
			\ddot{z}+z''''+B^2(Jz')'=0,& (x,t)\in (0,L)\times\R^+ \\
			z(0,t)=z'(0,t)=\phi(0,t)=0, &\\
			\phi'(L,t)=z''(L,t)=z'''(L,t)-B\phi(L,t)=0,& t\in  \R^+\\
			z(x,0)=z^0(x),\quad z_t(x,0)=z^1(x).& x\in[0,L]
		\end{array}\right.
	\end{eqnarray}
	Define $\mc{H} =\mathrm V \times \mathrm H= H^2_L(0,L) \times \Ltwo(0,L).$  Define also the  natural energy of the solutions of \eqref{compactclamped} as
	\begin{equation*}
		E(t):=\frac{1}{2}\int_{0}^{L}\left[|\dot{z}|^2+|z''|^2-B^2(Jz')z'\right]dx
	\end{equation*}
	where the term $-\frac{1}{2}\int_{0}^{L}B^2(Jz')z'~dx$ is non-negative due to \cref{lemma1}.
 This motivates the definition of the inner product on $\mc{H}:$
\begin{eqnarray}
\nonumber && \left<\left[ \begin{array}{l}
 u_1 \\
 u_2
 \end{array} \right], \left[ \begin{array}{l}
 v_1 \\
 v_2
 \end{array} \right]\right>_{\mc{H}}= \left<u_2, v_2\right>_{\mathrm H} + \left<u_1, v_1\right>_{\mathrm V}= \int_0^L \left\{m  u_2 { {\bar v}_2}+  (u_2)_{xx} (\bar v_2)_{xx}  - B^2 (J(u_1)_{x}) ({\bar u}_1)_x  \right\}~dx.
 \end{eqnarray} It is straightforward to show that $\frac{dE}{dt}=0, i.e.,  E(t)\equiv E(0)$ for all $t > 0.$

Define  the operator $\mc A: {\text{Dom}}(\mc A)\subset \mc H \to \mc H$
where $\mc A= \left[ {\begin{array}{*{20}c}
   0 & I \\
       -D_x^4-B^2 D_x J D_x&  0  \\
\end{array}} \right]$ and
\begin{eqnarray}
 \label{A-MM-newd}  \left.
\begin{array}{ll}
{\rm {Dom}}(\mc A) := \{ (z_1,z_2)\in \mc H, z_2\in H^2_L(0,L),~~(z_1)_{xxx}+B^2  J (z_1)_x \in H^1(0,L), ~~ (z_1)_{xx}(L)=0, &\\
\qquad\qquad \qquad\qquad   (z_1)_{xxx}(L)+ B^2  J (z_1)(L)=0 \}. &
 \end{array} \right.
\end{eqnarray}

Choosing the state $\Phi=[w,\dot w]^{\rm T},$ the   system (\ref{compactclamped})  can be put into the  state-space form
\begin{eqnarray}
\label{Semigroupp-mmm}
\dot \Phi = {\mc A}  \Phi, \quad \Phi(x,0) =  \Phi ^0.
\end{eqnarray}
The following well-posedness result is immediate from \cite{ozer_dynamic_2018}:
\begin{theorem}\label{w-pff}
The operator ${\mc A}$ in \eqref{Semigroupp-mmm}  is unitary on $\mc H,$ and it generates a $C_0$-semigroup of contractions on $\mc H.$ Letting $T>0,$ and for any $\Phi^0 \in \mathrm{H},$ $\Phi\in C[[0,T]; \mc H]$ and there exists a positive constant $c_1(T)$
such that (\ref{Semigroupp-mmm}) satisfies
      $\|\Phi (T) \|^2_{\mc{H}} \le c_1 (T)\|\Phi^0\|^2_{\mc H}.$
\end{theorem}

The first important result of this section, the  exact observability of  \eqref{compactclamped},  is provided in the following theorem.
	\begin{theorem}\label{cf-imp1}
		Let $T>2c$ and $c:=\max\left\{L,\frac{L^3}{\pi^2}\right\}.$ Then, the following observability inequality holds
	\begin{equation}\label{dmbb}
			\int_{0}^{T}|\dot{z}(L,t)|^2dt\geq \frac{2}{L}(T-2c)E(0).
		\end{equation}
	\end{theorem}
	\begin{proof} The proof is based on the multipliers technique, i.e. see \cite{komornik_exact_1997}.  First, multiply \eqref{compactclamped} by $xz'$ and integrate over $[0,L]$ with respect to $x$, and over $[0,T]$ with respect to $t$:
		\begin{align*}
			0=\int_{0}^{T}\int_{0}^{L}(\ddot{z}+z''''+B^2(Jz')')xz'dxdt=:I_1+I_2+I_3
		\end{align*}
		where
		\begin{align}
			I_1&=\int_{0}^{T}\int_{0}^{L}\ddot{z}xz'dxdt
=\left.\int_{0}^{L}\dot{z}xz'\right|_0^Tdx-\frac{L}{2}\int_{0}^{T}|\dot{z}(L,t)|^2dt+\frac{1}{2}\int_{0}^{T}\int_{0}^{L}|\dot{z}|^2 dxdt\\
					I_2&=\int_{0}^{T}\int_{0}^{L}z''''xz'dxdt
			=-LB^2\int_{0}^{T}(Jz'(L))z'(L)dt+\frac{3}{2}\int_{0}^{T}\int_{0}^{L}|z''|^2 dxdt,\\
\color{black} I_3&\color{black}=B^2\int_{0}^{T}\int_{0}^{L}(Jz')'xz'dxdt=\left.B^2\int_{0}^{T}(Jz')xz' \right|_0^Ldt-B^2\int_{0}^{T}\int_{0}^{L}(Jz')(z'+xz'')dxdt\nonumber\\
\label{dmb1} &=-B^2\int_{0}^{T}\int_{0}^{L}(Jz')xz''dxdt-B^2\int_{0}^{T}\int_{0}^{L}Jz'z'dxdt+LB^2\int_{0}^{T}(Jz'(L))z'(L)dt.
		\end{align}
		Let $(-CD_x^2+P)^{-1}z'=y$. Then, $z'=-Cy''+Py$. By the definition of the operator $(-CD_x^2+P)$, $y\in H^2_L$,  $y(0)=y'(L)=0$. Hence, the first integral in \eqref{dmb1} is rewritten as
		\begin{align*}
			-B^2\int_{0}^{T}\int_{0}^{L}(Jz')xz''dxdt&=B^2\int_{0}^{T}\int_{0}^{L}\left(\frac{z'}{C}-\frac{P}{C}(-CD_x^2+P)^{-1}z'\right)xz''dxdt\\
			&=B^2\int_{0}^{T}\int_{0}^{L}x(Cy'''y''-Py'y'')dxdt\\
			&=\frac{B^2C}{2}\int_{0}^{T}|y''(L)|^2dt-\frac{B^2C}{2}\int_{0}^{T}\int_{0}^{L}|y''|^2dxdt+\frac{B^2P}{2}\int_{0}^{T}\int_{0}^{L}|y'|^2 dxdt.
		\end{align*}
		The second integral in \eqref{dmb1} is also rewritten as
		\begin{align*}
			-\frac{B^2}{2}\int_{0}^{T}\int_{0}^{L}(Jz')z'dxdt&=\frac{B^2}{2}\int_{0}^{T}\int_{0}^{L}\left(\frac{z'}{C}-\frac{P}{C}(-CD_x^2+P)^{-1}z'\right)z'dxdt\\
			&=\frac{B^2}{2}\int_{0}^{T}\int_{0}^{L}-y''(-Cy''+Py)dxdt\\
			&=\frac{B^2C}{2}\int_{0}^{T}\int_{0}^{L}|y''|^2dxdt+\frac{B^2P}{2}\int_{0}^{T}\int_{0}^{L}|y'|^2dxdt.
		\end{align*}
		Therefore,
		\begin{eqnarray*}
			&&-B^2\int_{0}^{T}\int_{0}^{L}(Jz')xz''dxdt-\frac{B^2}{2}\int_{0}^{T}\int_{0}^{L}(Jz')z'dxdt =\frac{B^2C}{2}\int_{0}^{T}|y''(L)|^2dt+B^2P\int_{0}^{T}\int_{0}^{L}|y'|^2dxdt\geq 0.
		\end{eqnarray*}
		Next, the term $\left.\int_{0}^{L}\dot{z}xz'\right|_0^Tdx$ in $I_1$ is estimated. Since $c=\max\left\{L,\frac{L^3}{\pi^2}\right\},$
		\begin{align*}
			\left|\int_{0}^{L}\dot{z}xz'dx\right|&\leq \frac{L}{2}\int_{0}^{L}|z'|^2dx+\frac{L}{2}\int_{0}^{L}|\dot{z}|^2dx\leq \frac{L^3}{2\pi^2}\int_{0}^{L}|z''|^2dx+\frac{L}{2}\int_{0}^{L}|\dot{z}|^2dx \leq c E(0),
		\end{align*}
		and therefore, $\left.\int_{0}^{L}\dot{z}xz'\right|_0^Tdx  \geq -2c E(0).$ Finally,  \eqref{dmbb} follows by bringing all together $I_1+I_2+I_3=0:$
		\begin{align*}
			&\underbrace{\left.\int_{0}^{L}\dot{z}xz'\right|_0^Tdx}_{\geq -2cE(0)}-\frac{L}{2}\int_{0}^{T}|\dot{z}(L)|^2dt-\underbrace{B^2\int_{0}^{T}\int_{0}^{L}(Jz')xz''dxdt-\frac{B^2}{2}\int_{0}^{T}\int_{0}^{L}Jz'z'dxdt}_{\geq0}\\
			&+\underbrace{\int_{0}^{T}\int_{0}^{L}\frac{1}{2}|\dot{z}|^2+\frac{3}{2}|z''|^2-\frac{1}{2}B^2Jz'z'dxdt}_{\geq TE(0)}=0.\qedhere
		\end{align*}
 	\end{proof}
	
	\section{Semi-discretization in the Space Variable by the Order-reduced Finite-Differences}
	\label{clampeddiscretiation}
	In this section, the so-called ``Order-reduced Finite-Differences'' method, introduced in \cite{Guo2} for a single beam equation,  is used to approximate \eqref{eq1a}-\eqref{eq1d}. For the simplicity of the calculations, the length of beam is chosen $L=1$. The  choices of states
		$u(x,t):=\dot{z}(x,t),\quad v(x,t):=z'''(x,t)$ are needed
	so that the system \eqref{eq1a}-\eqref{eq1d} can be reformulated in the following first-order form (in time)
	\begin{equation}\label{orderreduced}
	\left\{	\begin{array}{ll}
			-\dot{v}+u'''=0& \\
			\dot{u}+v'-B\phi'=0&\\
			-C\phi''+P\phi=-Bv,\hspace{0.5cm}& (x,t)\in (0,1)\times\R^+\\
			u(0,t)=u'(0,t)=u''(1,t)=\phi(0,t)=v(1,t)-B\phi(1,t)=0,& t \in \mathbb{R}^+ \\
			v(x,0)=v^0(x),\quad u(x,0)=u^0(x),& x \in (0,L).
		\end{array}\right.
	\end{equation}
Let $N\in\N$ be given. Defining the mesh size $h:=\frac{1}{N+1}$, the following discretization of the interval $[0,L]$ is considered
	\begin{equation}\label{discret}
		0=x_0<x_1<...<x_{N-1}<x_N<x_{N+1}=L.
	\end{equation}
For this method, we also consider the middle points of each subinterval, and denote them by $\{x_{i+\frac{1}{2}}\}_{i=0}^{N+1},$ where i.e. $x_{i+\frac{1}{2}}=\frac{x_{i+1}+x_i}{2}.$ 	Define following finite average and difference operators
	\begin{align*}
		u_{i+\frac{1}{2}}&:=\frac{u_{i+1}+u_i}{2}, \quad \delta_xu_{i+\frac{1}{2}}:=\frac{u_{i+1}-u_i}{h}, \quad
		\delta_x^2u_{i}:=\frac{u_{i+1}-2u_i+u_{i-1}}{h^2}\\
		\delta_x^3u_{i+\frac{1}{2}}&:=\frac{u_{i+2}-3u_{i+1}+3u_i-u_{i-1}}{h^3}, \quad \delta_x^4u_{i}:=\frac{u_{i+2}-4u_{i+1}+6u_i-4u_{i-1}+u_{i-2} }{h^4}.
	\end{align*}
Note that considering the odd-number derivatives at the in-between nodes within the uniform discretization \eqref{discret} provides  higher-order approximations. Moreover, they can be represented in two separate ways by the preceding  and following nodes \cite{Guo1}.

The following application of difference operators to \eqref{orderreduced} yields a novel model reduction of  \eqref{orderreduced}. The discretization of the second equation of \eqref{orderreduced} is done as the following
	\begin{align}\label{orfd1a}
		&\dot{u}_{i+\frac{1}{2}}+\delta_xv_{i+\frac{1}{2}}-B\delta_x\phi_{i+\frac{1}{2}}=0.
	\end{align}
	Multiplying \eqref{orfd1a} with $\frac{h}{2}$ yields
	\begin{align}\label{orfd2a}
	&{\color{black}	\frac{v_{i+1}-v_i}{2}=-\frac{h}{2}\dot{u}_{i+\frac{1}{2}}+\frac{Bh}{2}\delta_x\phi_{i+\frac{1}{2}}.}
	\end{align}
	Observe that $\frac{v_{i+1}-v_i}{2}=v_{i+1}-v_{i+\frac{1}{2}}=v_{i+\frac{1}{2}}-v_i$ implies two representations of the left-hand side of \eqref{orfd2a} in terms of $v_{i+1}$ and $v_i$. First, substitute $\frac{v_{i+1}-v_i}{2}=v_{i+1}-v_{i+\frac{1}{2}}$ in \eqref{orfd2a} and shift indices by one, i.e. $i+1\to i,$
	\begin{align}
		v_{i+1}&=v_{i+\frac{1}{2}}-\frac{h}{2}\dot{u}_{i+\frac{1}{2}}+\frac{Bh}{2}\delta_x\phi_{i+\frac{1}{2}}, \quad
		v_{i}=v_{i-\frac{1}{2}}-\frac{h}{2}\dot{u}_{i-\frac{1}{2}}+\frac{Bh}{2}\delta_x\phi_{i-\frac{1}{2}}. \label{eq45}
	\end{align}
	and next, substitute$\frac{v_{i+1}-v_i}{2}=v_{i+\frac{1}{2}}-v_i$ in \eqref{orfd2a} to get
	\begin{equation}\label{eq46}
		v_{i}=v_{i+\frac{1}{2}}+\frac{h}{2}\dot{u}_{i+\frac{1}{2}}-\frac{Bh}{2}\delta_x\phi_{i+\frac{1}{2}},
	\end{equation}
	and then  subtract \eqref{eq45} from \eqref{eq46} to eliminate  $v_i$
	\begin{equation}\label{eq47}
		v_{i+\frac{1}{2}}-v_{i-\frac{1}{2}}+\frac{h}{2}\left(\dot{u}_{i+\frac{1}{2}}+\dot{u}_{i-\frac{1}{2}}\right)-\frac{Bh}{2}\left(\delta_x\phi_{i+\frac{1}{2}}+\delta_x\phi_{i-\frac{1}{2}}\right)=0.
	\end{equation}
	Finally, multiply \eqref{eq47} by $\frac{1}{h}$ and rewrite the equation in the $z-$form
	\begin{equation}\label{eq7}
		\delta_x^4z_i+\frac{1}{2}\left(\ddot{z}_{i+\frac{1}{2}}+\ddot{z}_{i-\frac{1}{2}}\right)-\frac{B}{2}\left(\delta_x\phi_{i+\frac{1}{2}}+\delta_x\phi_{i-\frac{1}{2}}\right)=0,
	\end{equation}
	which is equivalent to
	\begin{equation}\label{eq10}
		\delta_x^4z_i+\frac{1}{4}\left(\ddot{z}_{i+1}+2z_i+\ddot{z}_{i-1}\right)-\frac{B}{2h}\left(\phi_{i+1}-\phi_{i-1}\right)=0.
	\end{equation}
	Let $\tilde{y}_{i+\frac{1}{2}}:=\delta_x \phi_{i+\frac{1}{2}}$.  The third equation in \eqref{orderreduced} can be discretized as following
	\begin{align}\label{dmb5}
		-C\delta_x \tilde{y}_{i+\frac{1}{2}}+P\phi_{i+\frac{1}{2}}&=-Bv_{i+\frac{1}{2}}
	\end{align}
	or, by multiplying \eqref{dmb5} by $-\frac{h}{2},$
	\begin{align}\label{dmb2}
	\color{black}	C\left(\frac{\tilde{y}_{i+1}-\tilde{y}_i}{2}\right)& \color{black}=\frac{Ph}{2}\phi_{i+\frac{1}{2}}+\frac{Bh}{2}v_{i+\frac{1}{2}}.
	\end{align}
	Now, first, substitute $\frac{\tilde{y}_{i+1}-\tilde{y}_i}{2}=\tilde{y}_{i+\frac{1}{2}}-\tilde{y}_i$  in  \eqref{dmb2}
	\begin{equation}\label{eq31}
		C\tilde{y}_i=C\tilde{y}_{i+\frac{1}{2}}-\frac{Ph}{2}\phi_{i+\frac{1}{2}}-\frac{Bh}{2}v_{i+\frac{1}{2}},
	\end{equation}
	and next, substitute $\frac{\tilde{y}_{i+1}-\tilde{y}_i}{2}=\tilde{y}_{i+1}-\tilde{y}_{i+\frac{1}{2}}$ in  \eqref{dmb2}, and finally, shift indices by one i.e. $i+1\to i$  to get
	\begin{align}
		C\tilde{y}_{i+1}&=C\tilde{y}_{i+\frac{1}{2}}+\frac{Ph}{2}\phi_{i+\frac{1}{2}}+\frac{Bh}{2}v_{i+\frac{1}{2}}\\
		C\tilde{y}_{i}&=C\tilde{y}_{i-\frac{1}{2}}+\frac{Ph}{2}\phi_{i-\frac{1}{2}}+\frac{Bh}{2}v_{i-\frac{1}{2}}.\label{eq32}
	\end{align}
	Elimination of $\tilde{y}_i$ by the use of \eqref{eq31} and \eqref{eq32} yields,
	\begin{align*}
		&C(\tilde{y}_{i+\frac{1}{2}}-\tilde{y}_{i-\frac{1}{2}})-\frac{Ph}{2}(\phi_{i+\frac{1}{2}}+\phi_{i-\frac{1}{2}})-\frac{Bh}{2}(v_{i+\frac{1}{2}}+v_{i-\frac{1}{2}})=0\\
		&C(\delta_x\phi_{i+\frac{1}{2}}-\delta_x \phi_{i-\frac{1}{2}})-\frac{Ph}{2}(\phi_{i+\frac{1}{2}}+\phi_{i-\frac{1}{2}})-\frac{Bh}{2}(\delta_x^3z_{i+\frac{1}{2}}+\delta_x^3z_{i-\frac{1}{2}})=0\label{eq35}\\
		&C\left(\frac{\phi_{i+1}-2\phi_i+\phi_{i-1}}{h}\right)-\frac{Ph}{4}(\phi_{i+1}+2\phi_i+\phi_{i-1})-\frac{B}{2}\left(\frac{z_{i+2}-2z_{i+1}+2z_{i-1}-z_{i-2}}{h^2}\right)=0.
	\end{align*}
	Next, consider discretized boundary condition corresponding to $z'''(1)-B\phi(1)=0,$
	\begin{equation}\label{eq8}
		v_{N+1}-B\phi_{N+1}=0.
	\end{equation}
	By \eqref{eq47}, $v_{N+\frac{1}{2}}-\frac{h}{2}\dot{u}_{N+\frac{1}{2}}+\frac{Bh}{2}\delta_x\phi_{N+\frac{1}{2}}-B\phi_{N+1}=0,$
	and  the boundary condition $\delta_x\phi_{N+\frac{1}{2}}=0$, i.e. and $\phi_{N+1}=\phi_N,$
	\begin{equation}\label{eq9}
		-\delta_x^3z_{N+\frac{1}{2}}+\frac{h}{4}\left(\ddot{z}_{N+1}+\ddot{z}_{N}\right)+B\phi_{N}=0.
	\end{equation}
	Hence, the reduced-order discretization of \eqref{orderreduced} or \eqref{eq1a}-\eqref{eq1d} is as the following
	\begin{equation}\label{clampeddiscrete}
		\begin{cases}
			\delta_x^4z_i+\frac{1}{4}\left(\ddot{z}_{i+1}+2\ddot{z}_i+\ddot{z}_{i-1}\right)-\frac{B}{2h}\left(\phi_{i+1}-\phi_{i-1}\right)=0,&\\
			-C\left(\frac{\phi_{i+1}-2\phi_i+\phi_{i-1}}{h}\right)+\frac{Ph}{4}(\phi_{i+1}+2\phi_i+\phi_{i-1})+\frac{B}{2}\left(\frac{z_{i+2}-2z_{i+1}+2z_{i-1}-z_{i-2}}{h^2}\right)=0,& i=1,...,N,\\
			-\delta_x^3z_{N+\frac{1}{2}}+\frac{h}{4}\left(\ddot{z}_{N+1}+\ddot{z}_{N}\right)+B\phi_{N}=0,\\
			\phi_0=0, ~\phi_{N+1}=\phi_N, ~ z_0=z_{-1}=0, ~z_{N+2}=2z_{N+1}-z_N,\\
			z_i(0)=z^0(x_i),\quad \dot{z}_i(0)=z^1(x_i),& {\color{black}i=0,1,\ldots, N+1.}
		\end{cases}
	\end{equation}
	For convenience, let $y_i:=z_{i+1}-z_{i-1}$ so that the second equation can be rewritten as
	\begin{align}
	\label{dd7}	0&=-C\left(\frac{\phi_{i+1}-2\phi_i+\phi_{i-1}}{h}\right)+\frac{Ph}{4}(\phi_{i+1}+2\phi_i+\phi_{i-1})+\frac{B}{2}\left(\frac{y_{i+1}-2y_i+y_{i-1}}{h^2}\right).
	\end{align}
	The following results are needed for the rest of the section.
	\begin{lemma}\label{poincare}
		Let $\{u_i\}_{i=0}^{N+1}\in \R^{N+2}$ be the value of $u(x)$ at grids $\{x_i\}_{i=0}^{N+1}$ such that $u_0=0$ and $\delta_xu_{\frac{1}{2}}=0$. Then,
		\begin{equation}\nonumber
			{\color{black}\max\limits_{0\leq i\leq N+1}(u_i)^2\leq h\sum_{j=0}^{N}\left(\delta_x u_{j+\frac{1}{2}}\right)^2,\quad 	\max\limits_{0\leq j\leq N+1}\left(\delta_x u_{j+\frac{1}{2}}\right)^2\leq h\sum_{j=0}^{N}\left(\delta_x^2 u_{j}\right)^2.}
		\end{equation}
	\end{lemma}
	\begin{proof}
		Since $u_0=0$ and $\delta_xu_{\frac{1}{2}}=0$,
		\begin{equation*}
			u_i=(u_i-u_{i-1})+(u_{i-1}-u_{i-2})+...+(u_1-u_{0})=\sum_{j=0}^{i-1}(u_{j+1}-u_j)=h\sum_{j=0}^{i-1}\delta_x u_{j+\frac{1}{2}}, \quad i=1,...,N+1
		\end{equation*}
		and $\delta_x u_{i+\frac{1}{2}}=h\sum\limits_{j=1}^{i}\delta_x^2u_j.$
		By Cauchy-Schwarz's inequality,
		\begin{equation*}
\begin{array}{ll}
			(u_i)^2&=\left(h\sum\limits_{j=0}^{i-1}\delta_x u_{j+\frac{1}{2}} \right)^2\leq \left(h\sum_{j=0}^{i-1}1^2 \right)\left(h\sum\limits_{j=0}^{i-1}\left(\delta_x u_{j+\frac{1}{2}}\right)^2\right)
			=h\sum\limits_{j=0}^{N}\left(\delta_x u_{j+\frac{1}{2}}\right)^2\\
(\delta_x u_{j+\frac{1}{2}})^2&=\left(h\sum\limits_{j=1}^{i}\delta_x^2 u_{j} \right)^2\leq \left(h\sum\limits_{j=1}^{i}1^2 \right)\left(h\sum\limits_{j=1}^{i}\left(\delta_x^2 u_{j}\right)^2\right)
			=h\sum\limits_{j=1}^{N}\left(\delta_x^2 u_{j}\right)^2. \qedhere
\end{array}
		\end{equation*}

	\end{proof}
	\begin{corollary}\label{corallery}
		Assume the conditions in Lemma \ref{poincare}. Then, the solution $\{z_j\}_{j=0}^{N+1}$ of  \eqref{clampeddiscrete} satisfies \begin{equation*}
			\sum_{j=1}^{N}\left(z_{j}\right)^2\leq\sum_{j=1}^{N}\left(\delta_xz_{j+\frac{1}{2}}\right)^2\leq\sum_{j=1}^{N}\left(\delta_x^2z_{j}\right)^2.
		\end{equation*}
	\end{corollary}
	\begin{proof}
	By  $Nh\leq1,$ it is straightforward to show the second inequality
		\begin{align*}
			\sum_{j=1}^{N}\left(\delta_xz_{j+\frac{1}{2}}\right)^2&\leq N\max\limits_{0\leq j\leq N+1}\left(\delta_x z_{j+\frac{1}{2}}\right)^2\leq Nh\sum_{j=1}^{N}\left(\delta_x^2z_{j}\right)^2\leq\sum_{j=1}^{N}\left(\delta_x^2z_{j}\right)^2.
		\end{align*}
		The proof of the first inequality follows the similar steps.
	\end{proof}
	
	Assuming that $\delta_xz_{\frac{1}{2}}=0$ and $\delta_x^2z_N=0$ implies $z_1=0$, and $2z_{N}=z_{N+1}+z_{N-1}$. Therefore,
	\begin{align*}
		y_0&=z_1-z_{-1}=0, y_{N+1}=z_{N+2}-z_N=2(z_{N+1}-z_N)=z_{N+1}-z_{N-1}=y_N.
	\end{align*}
	For $i=1$ and $\phi_0=y_0=0,$
	\begin{align}\label{dd5}
		C\left(\frac{-\phi_{2}+2\phi_1}{h}\right)+\frac{Ph}{4}(\phi_{2}+2\phi_1)+\frac{B}{2}\left(\frac{y_{2}-2y_1}{h^2}\right)&=0,
	\end{align}
	and for $i=N,$ $\phi_{N+1}=\phi_N$ and $y_{N+1}=y_N,$
	\begin{align}
\label{dd6}		-C\left(\frac{-\phi_N+\phi_{N-1}}{h}\right)+\frac{Ph}{4}(3\phi_N+\phi_{N-1})+\frac{B}{2}\left(\frac{-y_N+y_{N-1}}{h^2}\right)&=0.
	\end{align}
	Thus, by  letting $\vec{\phi}=[\phi_1,...,\phi_{N}]^T$, $\vec{y}=[y_1,...,y_{N}]^T$, \eqref{dd7},\eqref{dd5},\eqref{dd6} can be rewritten in the following matrix form
	\begin{equation*}
		C\bm{A}_h\vec{\phi}+P\bm{M}\vec{\phi}-\frac{B}{2h}\bm{A}_h\vec{y}=0
	\end{equation*}
	where
	\begin{align*}
		\bm{A}_h=\frac{1}{h^2}\begin{bmatrix}
			2&-1&0&\dots&\dots&\dots&0\\
			-1&2&-1&0&\dots&\dots&0\\
			0&-1&2&-1&0&\dots&0\\
			&\ddots&\ddots&\ddots&\ddots&\ddots&\\
			0&\dots&0&1&-2&1&0\\
			0&\dots&\dots&0&-1&2&-1\\
			0&\dots&\dots&\dots&0&-1&1\\
		\end{bmatrix},~~\bm{M}=\frac{1}{4}\begin{bmatrix}
			2&1&0&\dots&\dots&\dots&0\\
			1&2&1&0&\dots&\dots&0\\
			0&1&2&1&0&\dots&0\\
			&\ddots&\ddots&\ddots&\ddots&\ddots&\\
			0&\dots&0&1&2&1&0\\
			0&\dots&\dots&0&1&2&1\\
			0&\dots&\dots&\dots&0&1&3\\
		\end{bmatrix}.
	\end{align*}
	This immediately leads to a better representation of $\vec{\phi}=\frac{B}{2h}\left(C\bm{A}_h+P\bm{M}\right)^{-1}\bm{A}_h\vec{y}.$
	The following result is crucial to show that the matrices $\bm{A}_h$ and $\bm{M}$ have the same eigenvectors.
	\begin{lemma}\label{matrices}
		The matrices $\bm{A}_h$ and $\bm{M}$ are positive-definite, symmetric matrices. Moreover, eigenvectors of both $\bm{M}$ and $\bm{A}_h$ are
		\begin{align*}
			\varphi^k&=\begin{bmatrix}
				\varphi_{k,1}&\cdots&\varphi_{k,N}
			\end{bmatrix}^T, \quad \varphi_{k,j}=\sin\left(\frac{(2j-1)k\pi}{2N+1}\right), \quad  k,j=1,\ldots, N.
		\end{align*}
	\end{lemma}
	\begin{proof}
		For $u=\begin{bmatrix}
			u_1&u_2&\cdots&u_N
		\end{bmatrix}^T$, solving  the eigenvalue problem $\bm{A}_hu=\lambda u$ is equivalent to solving the  system of homogeneous difference equations
		\begin{equation*}
			\begin{cases}
				-u_{k+1}+(2-h^2\lambda)u_{k}-u_{k-1}=0, & k=1,2,...,N\\
				u_0=0, \quad u_{N+1}=u_N.
			\end{cases}
		\end{equation*}
		Let $z^k:=u_k$ be a solution. Then,
			{\color{black}$z^{k-1}\left(-1+(2-h^2\lambda)z-z^2\right)=0.$}
		Since $z^{k-1}\neq 0$, $-1+(2-h^2\lambda)z-z^2=0$. Thus, there exist arbitrary constants $c_1,c_2$ such that the general solution of the difference equation above can be expressed as
		\begin{equation*}
			u_k=c_1z^k+c_2z^{-k}.
		\end{equation*}
		By the boundary condition $u_0=0,$  $-c_1=c_2$. By the other boundary condition $u_{N+1}=u_N,$
		\begin{align*}
			c_1\left(z^{N+1}-z^{-(N+1)}\right)&=c_1\left(z^N-z^{-N}\right)\\
			z^{N+1}-z^N&=z^{-(N+1)}-z^{-N}\\
			z^{2N+1}&=-1=e^{i(2j-1)\pi},\quad j=1,2,...,N.
		\end{align*}
		Therefore, $z_j=e^{\frac{i(2j-1)\pi}{2N+1}}, j=1,2,...,N,$
		and the eigenvectors of ${\bf A_h}$ are
		\begin{align*}
			\varphi^k&=\begin{bmatrix}
				\varphi_{k,1}&\cdots&\varphi_{k,N}
			\end{bmatrix}^T,\quad
			\varphi_{k,j}=\sin\left(\frac{(2j-1)k\pi}{2N+1}\right), ~~k,j=1,2,...,N.
		\end{align*}
		Next, the eigenvalues are solved by $z_1+z_2=2-h\lambda_k$,
		\begin{align*}
			2\left(1-2\sin^2\left( \frac{(2k-1)\pi}{4N+2}\right)\right)&=2-h^2\lambda_k\\
			\lambda_k&=\frac{4}{h^2}\sin^2\left( \frac{(2k-1)\pi}{4N+2}\right), \quad k=1,2,...,N.
		\end{align*}
		Observe that the matrices $M$ and $A_h$ satisfy 	$h^2\bm{A_h}=4I-4\bm{M},$
		which implies that the eigenvectors of $\bm{M}$ and $\bm{A_h}$ are the same, and the eigenvalues of $\bm{M}$ and $\bm{M^{-1}A_h}$, respectively, are explicitly written as
		\begin{equation*}
			\tilde{\lambda}_k=1-\sin^2\left( \frac{(2k-1)\pi}{4N+2}\right), \quad
			\mu_k=\frac{4\sin^2\left( \frac{(2k-1)\pi}{4N+2}\right)}{h^2-h^2\sin^2\left( \frac{(2k-1)\pi}{4N+2}\right)}, \quad k=1,2,...,N.
		\end{equation*}
		Finally, the matrices  $M$ and $A_h$ are positive definite since their eigenvalues are all positive and real.
	\end{proof}
	\section{Conservation of the Discretized Energy}
	Define the discrete energy of the system as
	\begin{equation}\label{discreteE}
		E_h(t):=\frac{h}{2}\sum_{j=0}^{N}\left(\dot{z}_{j+\frac{1}{2}}\right)^2+\frac{h}{2}\sum_{j=0}^{N}\left(\delta_x^2z_j\right)^2+\frac{B}{4}\sum_{j=0}^{N}\phi_j(z_{j+1}-z_{j-1}).
	\end{equation}
	The following results  are crucial for the rest of the section.
	\begin{lemma}\label{joperator}
		Define $\bm{J}_h:=\frac{1}{C}\bm{I}-\frac{P}{C}\left(C\bm{A}_h+P\bm{M}\right)^{-1}\bm{M}$. Then, $\bm{J}_h$ is a  self-adjoint operator. Moreover, for all $u\in \text{Dom}(\bm{J}_h),$
		\begin{equation*}
			\bm{J}_h=\left(\bm{C}\bm{A}_h+P\bm{M}\right)^{-1}\bm{A}_h.
		\end{equation*}
	\end{lemma}

\begin{proof}
		Since $\bm{A}_h$ and $\bm{M}$ are symmetric and positive definite real matrices with the same eigenvectors, $\bm{J}_h$ is self-adjoint, i.e. $\bm{J}_h^T=\bm{J}_h$. Let $\frac{1}{C}u-\frac{P}{C}\left(C\bm{A}_h+P\bm{M}\right)^{-1}\bm{M}u=v.$ Then, $-\bm{M}u=\frac{C^2}{P}\bm{A}_hv+C\bm{M}v-	\frac{C}{P}\bm{A}_hu-\bm{M}u$ and this trivially leads to $	\left(\bm{C}\bm{A}_h+P\bm{M}\right)^{-1}\bm{A}_hu=v.$
	\end{proof}
	Note that the second equation of \eqref{clampeddiscrete} and \cref{joperator}, the following relation between $\vec \phi$ and $\vec y$ is immediate:
	\begin{align}\label{dd1}
		\vec{\phi}=&\frac{B}{2h}\left(C\bm{A}_h+P\bm{M}\right)^{-1}\bm{A}_h\vec{y}=\frac{B}{2hC}\left(I-\left(\frac{C}{P}\bm{M}^{-1}\bm{A}_h+I\right)^{-1}\right)\vec{y}.
	\end{align}
	Now, the following intermediate variables are needed. Letting $\vec{u}:=\left(\frac{C}{P}\bm{M}^{-1}\bm{A}_h+I\right)^{-1}\vec{y}$
	implies that $	\vec{y}=\frac{C}{P}\bm{M}^{-1}\bm{A}_h\vec{u}+\vec{u}.$
	Substituting  $\vec y$ in \eqref{dd1} leads to
	\begin{align*}
		\vec{\phi}=&\frac{B}{2hC}\left(\vec{y}-\vec{u}\right)=\frac{B}{2Ph}\bm{M}^{-1}\bm{A}_h\vec{u}.
	\end{align*}
	Next, define $\vec{k}:=\bm{M}^{-1}\vec{u}$. Then, $
		\vec{y}=\frac{C}{P}\bm{A}_h\vec{k}+\bm{M}\vec{k},$	and by the relation $I-\frac{h^2}{4}\bm{A}_h=\bm{M},$ better forms of $\vec \phi$ and $\vec y$ are obtained
	\begin{equation}\label{kdef}
		\vec{\phi}=\frac{B}{2Ph}\bm{A}_h\vec{k},\quad\text{ and }\quad \vec{y}=\left(\frac{C}{P}-\frac{h^2}{4}\right)\bm{A}_h\vec{k}+\vec{k}.
	\end{equation}
The boundary conditions $\phi_0=0$ and $\phi_{N+1}=\phi_N$ imply that $\left(\bm{A}_h\vec{k} \right)_0=0$ and \linebreak $\left(\bm{A}_h\vec{k}\right)_{N+1}=\left(\bm{A}_h\vec{k}\right)_N$, respectively. The boundary conditions $y_0=0$ and $y_{N+1}=y_N$ also lead to
$y_0=\left(\frac{C}{P}-\frac{h^2}{4}\right)\left(\bm{A}_h\vec{k}\right)_0+k_0=k_0=0,$ and
\begin{eqnarray}
\begin{array}{ll}
	0=y_{N+1}-y_N=\left(\frac{C}{P}-\frac{h^2}{4}\right)\left(\bm{A}_h\vec{k}\right)_{N+1}+k_{N+1}-\left(\frac{C}{P}-\frac{h^2}{4}\right)\left(\bm{A}_h\vec{k}\right)_{N}-k_{N}=k_{N+1}-k_N,
\end{array}
\end{eqnarray}
and therefore, a new set of  boundary conditions on the auxiliary variable $\vec{k}$ is obtained as the following
\begin{equation}\label{boundary}
	k_0=\delta_xk_{N+\frac{1}{2}}=0.
\end{equation}
	\begin{lemma}\label{adjoint}
	Assume \eqref{kdef}-\eqref{boundary}.	Let $\{\phi_j, z_j\}_{j=0}^{N+1}$ be a solution of \eqref{clampeddiscrete}. Then, the following equality holds true
		\begin{equation}
			\sum_{j=0}^{N}\dot\phi_{j}({z}_{j+1}-{z}_{j-1})=\sum_{j=0}^{N}\phi_{j}(\dot{z}_{j+1}-\dot{z}_{j-1}).
		\end{equation}
	\end{lemma}
	\begin{proof}
		By $y_j=z_{j+1}-z_{j-1}$, \cref{matrices} and \eqref{kdef},
		\begin{align*}
			\sum_{j=0}^{N}\dot\phi_{j}({z}_{j+1}-{z}_{j-1})=&\sum_{j=0}^{N}\frac{B}{2Ph}\left(\bm{A}_h\dot{\vec{k}}\right)_j\left(\left(\frac{C}{P}-\frac{h^2}{4}\right)\left(\bm{A}_h\vec{k}\right)_j+\vec{k}_j\right)
			=\sum_{j=0}^{N}\phi_{j}(\dot{z}_{j+1}-\dot{z}_{j-1}). \qedhere
		\end{align*} 
	\end{proof}
	\begin{lemma}\label{conservedclamped}
		The discrete energy \eqref{discreteE} of the discrete system \eqref{clampeddiscrete} is conservative, i.e. $\frac{dE_h(t)}{dt}=0.$
		Hence, $E_h(t)=E_h(0)$ for all $t> 0$.
	\end{lemma}
	\begin{proof}
		First, multiply both sides of first equation of the system by $h\dot{z}_j$ and sum up $j$ from $1$ to $N$ to obtain
		\begin{equation}
			\frac{h}{2}\sum_{j=1}^{N}\left(\ddot{z}_{j+\frac{1}{2}}+\ddot{z}_{j-\frac{1}{2}}\right)\dot{z}_j+h\sum_{j=1}^{N}\delta_x^4z_j\dot{z}_j-\frac{B}{2}\sum_{j=1}^{N}\left(\phi_{i+1}-\phi_{i-1}\right)\dot{z}_j=0.
		\end{equation}
		Now define
		\begin{align*}
			I_1:=\frac{h}{2}\sum_{j=1}^{N}\left(\ddot{z}_{j+\frac{1}{2}}+\ddot{z}_{j-\frac{1}{2}}\right)\dot{z}_j,\quad
			I_2:=h\sum_{j=1}^{N}\delta_x^4z_j\dot{z}_j, \quad
			I_3:=-\frac{B}{2}\sum_{j=1}^{N}\left(\phi_{j+1}-\phi_{j-1}\right)\dot{z}_j
		\end{align*}
		such that $I_1+I_2+I_3=0$. By $\dot{z}_0=0$ in \eqref{clampeddiscrete},
		\begin{align*}
			I_1=\frac{h}{2}\sum_{j=1}^{N}\left(\ddot{z}_{j+\frac{1}{2}}+\ddot{z}_{j-\frac{1}{2}}\right)\dot{z}_j&=\frac{h}{2}\sum_{j=1}^{N}\ddot{z}_{j+\frac{1}{2}}\dot{z}_j+\frac{h}{2}\sum_{j=0}^{N-1}\ddot{z}_{j+\frac{1}{2}}\dot{z}_{j+1}\\
			&=\frac{h}{2}\sum_{j=1}^{N}\ddot{z}_{j+\frac{1}{2}}\dot{z}_j+\underbrace{\frac{h}{2}\ddot{z}_{\frac{1}{2}}\dot{z}_0}_{=0}+\frac{h}{2}\sum_{j=0}^{N-1}\ddot{z}_{j+\frac{1}{2}}\dot{z}_{j+1}+\underbrace{\frac{h}{2}\ddot{z}_{N+\frac{1}{2}}\dot{z}_{N+1}-\frac{h}{2}\ddot{z}_{N+\frac{1}{2}}\dot{z}_{N+1}}_{=0}\\
			&=h\sum_{j=0}^{N}\ddot{z}_{j+\frac{1}{2}}\dot{z}_{j+\frac{1}{2}}-\frac{h}{2}\ddot{z}_{N+\frac{1}{2}}\dot{z}_{N+1}.
		\end{align*}
		By $\dot{z}_0=0$, $\delta_x\dot{z}_{-\frac{1}{2}}=0$, and $\delta_x^2z_{N+1}=0$ in \eqref{clampeddiscrete},
		\begin{align*}
			I_2=h\sum_{j=1}^{N}\delta_x^4z_j\dot{z}_j&=\sum_{j=1}^{N}\delta_x^3z_{j+\frac{1}{2}}\dot{z}_j-\sum_{j=1}^{N}\delta_x^3z_{j-\frac{1}{2}}\dot{z}_j\\
			&=\sum_{j=1}^{N}\delta_x^3z_{j+\frac{1}{2}}\dot{z}_j-\sum_{j=0}^{N-1}\delta_x^3z_{j+\frac{1}{2}}\dot{z}_{j+1}\\
			&=\sum_{j=1}^{N}\delta_x^3z_{j+\frac{1}{2}}\dot{z}_j+\underbrace{\delta_x^3z_{\frac{1}{2}}\dot{z}_0}_{=0}-\sum_{j=0}^{N-1}\delta_x^3z_{j+\frac{1}{2}}\dot{z}_{j+1}+\underbrace{ \delta_x^3z_{N+\frac{1}{2}}\dot{z}_{N+1}-\delta_x^3z_{N+\frac{1}{2}}\dot{z}_{N+1}}_{=0}\\
			&=-\sum_{j=0}^{N}\left(\delta_x^2z_{j+1}-\delta_x^2z_{j}\right)\delta_x\dot{z}_{j+\frac{1}{2}}+\delta_x^3z_{N+\frac{1}{2}}\dot{z}_{N+1}\\
			&=-\sum_{j=1}^{N+1}\delta_x^2z_{j}\delta_x\dot{z}_{j-\frac{1}{2}}+\sum_{j=0}^{N}\delta_x^2z_{j}\delta_x\dot{z}_{j+\frac{1}{2}}+\delta_x^3z_{N+\frac{1}{2}}\dot{z}_{N+1}\\
			&=-\sum_{j=1}^{N+1}\delta_x^2z_{j}\delta_x\dot{z}_{j-\frac{1}{2}}-\underbrace{\delta_x^2z_{0}\delta_x\dot{z}_{-\frac{1}{2}}}_{=0}+\underbrace{\delta_x^2z_{N+1}\delta_x\dot{z}_{N+\frac{1}{2}}}_{=0}+\sum_{j=0}^{N}\delta_x^2z_{j}\delta_x\dot{z}_{j+\frac{1}{2}}+\delta_x^3z_{N+\frac{1}{2}}\dot{z}_{N+1}\\
			&=\sum_{j=0}^{N}\delta_x^2z_{j}\left(\delta_x\dot{z}_{j+\frac{1}{2}}-\delta_x\dot{z}_{j-\frac{1}{2}}\right)+\delta_x^3z_{N+\frac{1}{2}}\dot{z}_{N+1}\\
			&=h\sum_{j=0}^{N}\delta_x^2z_{j}\delta_x^2\dot{z}_{j}+\delta_x^3z_{N+\frac{1}{2}}\dot{z}_{N+1}.
		\end{align*}
		By all the boundary conditions in \eqref{clampeddiscrete},
		\begin{align*}
			I_3=&-\frac{B}{2}\sum_{j=1}^{N}\left(\phi_{j+1}-\phi_{j-1}\right)\dot{z}_j\\
			=&-\frac{B}{2}\sum_{j=2}^{N+1}\phi_{j}\dot{z}_{j-1}+\frac{B}{2}\sum_{j=0}^{N-1}\phi_{j}\dot{z}_{j+1}\\
			=&-\frac{B}{2}\sum_{j=0}^{N+1}\phi_{j}\dot{z}_{j-1}+\frac{B}{2}\sum_{j=0}^{N+1}\phi_{j}\dot{z}_{j+1}-\frac{3B}{2}\phi_N\dot{z}_{N+1}+\frac{B}{2}\phi_{N}\dot{z}_{N}\\
			=&\frac{B}{2}\sum_{j=0}^{N-1}\phi_{j}(\dot{z}_{j+1}-\dot{z}_{j-1})-\frac{B}{2}\phi_N(-\dot{z}_{N+2}+2\dot{z}_{N+1}+\dot{z}_{N-1})\\
			=&\frac{B}{2}\sum_{j=0}^{N-1}\phi_{j}(\dot{z}_{j+1}-\dot{z}_{j-1})+\frac{B}{2}\phi_N(\dot{z}_{N+2}-\dot{z}_{N-1})-B\phi_N\dot{z}_{N+1}\\
			=&\frac{B}{2}\sum_{j=0}^{N}\phi_{j}(\dot{z}_{j+1}-\dot{z}_{j-1})-B\phi_N\dot{z}_{N+1}.
		\end{align*}
		Moreover, the boundary equation in \eqref{clampeddiscrete} is multiplied by $\dot{z}_{N+1}$  to obtain
		\begin{equation*}
			\delta_x^3z_{N+\frac{1}{2}}\dot{z}_{N+1}-\frac{h}{2}\ddot{z}_{N+\frac{1}{2}}\dot{z}_{N+1}-B\phi_{N}\dot{z}_{N+1}=0.
		\end{equation*}
		Taking derivative of the energy with respect to $t$ and using \cref{adjoint} leads to
		\begin{align*}
			\frac{dE_h(t)}{dt}&=h\sum_{j=0}^{N}\ddot{z}_{j+\frac{1}{2}}\dot{z}_{j+\frac{1}{2}}+h\sum_{j=0}^{N}\delta_x^2z_{j}\delta_x^2\dot{z}_{j}+\frac{B}{4}\sum_{j=0}^{N}\dot\phi_{j}({z}_{j+1}-{z}_{j-1})+\frac{B}{4}\sum_{j=0}^{N}\phi_{j}(\dot{z}_{j+1}-\dot{z}_{j-1})\\
			&=h\sum_{j=0}^{N}\ddot{z}_{j+\frac{1}{2}}\dot{z}_{j+\frac{1}{2}}+h\sum_{j=0}^{N}\delta_x^2z_{j}\delta_x^2\dot{z}_{j}+\frac{B}{2}\sum_{j=0}^{N}\phi_{j}(\dot{z}_{j+1}-\dot{z}_{j-1})
		\end{align*}
		Hence,
		\begin{align*}
			0=I_1+I_2+I_3=&h\sum_{j=0}^{N}\ddot{z}_{j+\frac{1}{2}}\dot{z}_{j+\frac{1}{2}}-\frac{h}{2}\ddot{z}_{N+\frac{1}{2}}\dot{z}_{N+1}+h\sum_{j=0}^{N}\delta_x^2z_{j}\delta_x^2\dot{z}_{j}+\delta_x^3z_{N+\frac{1}{2}}\dot{z}_{N+1}\\
			&+\frac{B}{2}\sum_{j=0}^{N}\phi_{j}(\dot{z}_{j+1}-\dot{z}_{j-1})-B\phi_N\dot{z}_{N+1}\\
			&=\frac{dE_h(t)}{dt}+\delta_x^3z_{N+\frac{1}{2}}\dot{z}_{N+1}-\frac{h}{2}\ddot{z}_{N+\frac{1}{2}}\dot{z}_{N+1}-B\phi_{N}\dot{z}_{N+1}=\frac{dE_h(t)}{dt}.\qedhere
		\end{align*}
	\end{proof}
	
	\section{Uniform Observability by the Discrete Multipliers  as $h\to 0$}
	The following theorem is the main result  of this section.
	\begin{theorem}\label{mainthm}
Assume \eqref{kdef}-\eqref{boundary}, and
\begin{equation}\label{rmk}
			\left(\frac{C}{P}-\frac{h^2}{4}\right)-\frac{5}{2}\frac{B^2}{P}\geq 0.
		\end{equation} 	
Letting $\{\phi_j, z_j\}_{j=0}^{N+1}$ be a solution of \eqref{clampeddiscrete} and $T>6, $	the following observability inequality holds true
		\begin{equation}\label{MAIN}
		\int_0^T (\dot z_{N+1})^2~dt  \geq (T-6) E_h(t).
		\end{equation}
	\end{theorem}
	
	\begin{remark} Note that  $h$ in \eqref{rmk} is the mesh parameter of the discretization, and it satisfies $0<h<<1.$ The condition \eqref{rmk} simply reduces to $\frac{B^2}{C}<<\frac{2}{5},$ which is equivalent to
	$
		G_2<< \frac{2h_2(D_1h_1^3+D_3h_3^3)}{5H^2}.$
		In terms of the material parameters, one can interpret this condition as the large shear due to middle layer (small $G_2$). Indeed, this  is also desired in practical applications for more shear damping due to the middle layer.
	\end{remark}
	\begin{proof}
		In the proof, two discrete multipliers $hj\frac{z_{j+1}-z_{j-1}}{2}$ and $\frac{h}{4}z_j$ are separately used for only the first equation of \eqref{clampeddiscrete}. Correspondingly, two  multipliers $\delta_xz_{N+\frac{1}{2}}$ and $\frac{1}{4}z_{N+1}$ are separately used for third (boundary) equation in  \eqref{clampeddiscrete}. For the sake of convenience, first, multiply the first equation in \eqref{clampeddiscrete} by $hj\frac{z_{j+1}-z_{j-1}}{2}$, take the sum for $j=1,2,\ldots, N,$  and integrate by parts with respect to $t$ over $[0,T]$
		\begin{align}\label{denk1}
			\nonumber 0=&\frac{h}{2}\sum_{j=1}^{N}\int_{0}^{T}\left(\ddot{z}_{j+\frac{1}{2}}+\ddot{z}_{j-\frac{1}{2}}\right)j\frac{z_{j+1}-z_{j-1}}{2}dt+h\sum_{j=1}^{N}\int_{0}^{T}\delta_x^4z_jj\frac{z_{j+1}-z_{j-1}}{2}dt\\
			\nonumber &\qquad-\frac{B}{4}\sum_{j=1}^{N}\int_{0}^{T}(\phi_{j+1}-\phi_{j-1})j(z_{j+1}-z_{j-1})dt\\
\nonumber
		=&\frac{h}{2}\sum_{j=1}^{N}\left.\left(\dot{z}_{j+\frac{1}{2}}+\dot{z}_{j-\frac{1}{2}}\right)j\frac{z_{j+1}-z_{j-1}}{2}\right|_0^T-\frac{h}{2}\sum_{j=1}^{N}\int_{0}^{T}\left(\dot{z}_{j+\frac{1}{2}}+\dot{z}_{j-\frac{1}{2}}\right)j\frac{\dot{z}_{j+1}-\dot{z}_{j-1}}{2}dt\\
			&+h\sum_{j=1}^{N}\int_{0}^{T}\delta_x^4z_jj\frac{z_{j+1}-z_{j-1}}{2}dt-\frac{B}{4}\sum_{j=1}^{N}\int_{0}^{T}(\phi_{j+1}-\phi_{j-1})jy_j.
		\end{align}
	Now, multiply the third equation in \eqref{clampeddiscrete} by $\delta_xz_{N+\frac{1}{2}}$ and integrate with respect to $t$ over $[0,T]$
		\begin{align}
			\nonumber 0=&\frac{h}{2}\int_{0}^{T}\ddot{z}_{N+\frac{1}{2}}\delta_xz_{N+\frac{1}{2}}dt-\int_{0}^{T}\delta_x^3z_{N+\frac{1}{2}}\delta_xz_{N+\frac{1}{2}}dt +B\int_{0}^{T}\phi_N\delta_xz_{N+\frac{1}{2}}dt\\
			 \label{denk2} =&\left.\frac{h}{2}\dot{z}_{N+\frac{1}{2}}\delta_xz_{N+\frac{1}{2}}\right|_0^T-\frac{h}{2}\int_{0}^{T}\dot{z}_{N+\frac{1}{2}}\delta_x\dot{z}_{N+\frac{1}{2}}dt-\int_{0}^{T}\delta_x^3z_{N+\frac{1}{2}}\delta_xz_{N+\frac{1}{2}}dt+B\int_{0}^{T}\phi_N\delta_xz_{N+\frac{1}{2}}dt.
		\end{align}
	Next, define
		\begin{align}
\nonumber 			L_{h,1}(t)&:=\frac{h}{2}\sum_{j=1}^{N}\left(\dot{z}_{j+\frac{1}{2}}+\dot{z}_{j-\frac{1}{2}}\right)j\frac{z_{j+1}-z_{j-1}}{2}+\frac{h}{2}\dot{z}_{N+\frac{1}{2}}\delta_xz_{N+\frac{1}{2}},\\
		\nonumber 	I_1&:=-\frac{h}{2}\sum_{j=1}^{N}\int_{0}^{T}\left(\dot{z}_{j+\frac{1}{2}}+\dot{z}_{j-\frac{1}{2}}\right)j\frac{\dot{z}_{j+1}-\dot{z}_{j-1}}{2}dt-\frac{h}{2}\int_{0}^{T}\dot{z}_{N+\frac{1}{2}}\delta_x\dot{z}_{N+\frac{1}{2}}dt,\\
			\nonumber I_2&:=h\sum_{j=1}^{N}\int_{0}^{T}\delta_x^4z_jj\frac{z_{j+1}-z_{j-1}}{2}dt-\int_{0}^{T}\delta_x^3z_{N+\frac{1}{2}}\delta_xz_{N+\frac{1}{2}}dt,\\
			\label{denk3} I_3&:=-\frac{B}{4}\sum_{j=1}^{N}\int_{0}^{T}(\phi_{j+1}-\phi_{j-1})jy_jdt+B\int_{0}^{T}\phi_N\delta_xz_{N+\frac{1}{2}}dt,
		\end{align}
		such that the sum of \eqref{denk1} and \eqref{denk2} is
		\begin{equation}
\label{denk5}
			L_{h,1}(t)|_0^T+I_1+I_2+I_3=0.
		\end{equation}
		It is straightforward by the Young's inequality that
		\begin{equation}\label{denk6} L_{h,1}(t)\leq\frac{h}{8}\sum_{j=1}^{N}\left(\dot{z}_{j+\frac{1}{2}}+\dot{z}_{j-\frac{1}{2}}\right)^2+\frac{h}{8}\sum_{j=1}^{N}j^2\left(z_{j+1}-z_{j-1}\right)^2+\frac{h}{4}\left(\dot{z}_{N+\frac{1}{2}}\right)^2+\frac{h}{4}\left(\delta_xz_{N+\frac{1}{2}}\right)^2.
		\end{equation}
		Recalling that $z_{j+1}-z_{j-1}=h(\delta_xz_{j+\frac{1}{2}}+\delta_xz_{j-\frac{1}{2}})$ and $Nh<1$, the mean inequality $(a+b)^2\leq 2(a^2+b^2)$ implies  that $\left(\dot{z}_{j+\frac{1}{2}}+\dot{z}_{j-\frac{1}{2}}\right)^2\leq 2\left(\dot{z}_{j+\frac{1}{2}}\right)^2+2\left(\dot{z}_{j-\frac{1}{2}}\right)^2$. Therefore, the inequality \eqref{denk6} is majorized further
		\begin{align*}
			L_{h,1}(t)&\leq \frac{h}{4}\sum_{j=1}^{N}\left(\dot{z}_{j+\frac{1}{2}}\right)^2+\frac{h}{4}\sum_{j=1}^{N}\left(\dot{z}_{j-\frac{1}{2}}\right)^2+\frac{h}{8}\sum_{j=1}^{N}j^2h^2\left(\delta_xz_{j+\frac{1}{2}}+\delta_xz_{j-\frac{1}{2}}\right)^2+\frac{h}{4}\left(\dot{z}_{N+\frac{1}{2}}\right)^2+\frac{h}{4}\left(\delta_xz_{N+\frac{1}{2}}\right)^2\\
			\leq&\frac{h}{4}\sum_{j=1}^{N}\left(\dot{z}_{j+\frac{1}{2}}\right)^2+\left(\frac{h}{4}\sum_{j=0}^{N-1}\left(\dot{z}_{j+\frac{1}{2}}\right)^2+\frac{h}{4}\left(\dot{z}_{N+\frac{1}{2}}\right)^2\right)+\frac{h(N^2h^2)}{8}\sum_{j=1}^{N}\left(\delta_xz_{j+\frac{1}{2}}+\delta_xz_{j-\frac{1}{2}}\right)^2+\frac{h}{4}\left(\delta_xz_{N+\frac{1}{2}}\right)^2\\
			&\leq\frac{h}{2}\sum_{j=0}^{N}\left(\dot{z}_{j+\frac{1}{2}}\right)^2+\frac{h}{4}\sum_{j=1}^{N}\left(\delta_xz_{j+\frac{1}{2}}\right)^2+\frac{h}{4}\sum_{j=1}^{N}\left(\delta_xz_{j-\frac{1}{2}}\right)^2+\frac{h}{4}\left(\delta_xz_{N+\frac{1}{2}}\right)^2\\
			&\leq\frac{h}{2}\sum_{j=0}^{N}\left(\dot{z}_{j+\frac{1}{2}}\right)^2+\frac{h}{4}\sum_{j=1}^{N}\left(\delta_xz_{j+\frac{1}{2}}\right)^2+\left(\frac{h}{4}\sum_{j=0}^{N-1}\left(\delta_xz_{j+\frac{1}{2}}\right)^2+\frac{h}{4}\left(\delta_xz_{N+\frac{1}{2}}\right)^2\right)\\
			&\leq\frac{h}{2}\sum_{j=0}^{N}\left(\dot{z}_{j+\frac{1}{2}}\right)^2+\frac{h}{2}\sum_{j=1}^{N}\left(\delta_x^2z_{j}\right)^2\leq E_h(t).
		\end{align*}
	By the conservation of energy, \cref{conservedclamped}, $E_h(t)=E_h(0)$ for all $t\geq 0$, and thus,
		\begin{equation}\label{eq3}
			-2E_h(0)\leq L_{h,1}(t)|_0^T\leq 2E_h(0).
		\end{equation}
		Now, we turn our attention to $I_1$ in \eqref{denk3}. Since $\frac{\dot{z}_{j+1}-\dot{z}_{j-1}}{2}=\dot{z}_{j+\frac{1}{2}}-\dot{z}_{j-\frac{1}{2}},$
		\begin{align}
		\nonumber 	I_1&=-\frac{h}{2}\sum_{j=1}^{N}\int_{0}^{T}\left(\dot{z}_{j+\frac{1}{2}}+\dot{z}_{j-\frac{1}{2}}\right)j\frac{\dot{z}_{j+1}-\dot{z}_{j-1}}{2}dt-\frac{h}{2}\int_{0}^{T}\dot{z}_{N+\frac{1}{2}}\delta_x\dot{z}_{N+\frac{1}{2}}dt\\
	\nonumber 		&=-\frac{h}{2}\sum_{j=1}^{N}\int_{0}^{T}\left(\dot{z}_{j+\frac{1}{2}}\right)^2jdt+\frac{h}{2}\sum_{j=1}^{N}\int_{0}^{T}\left(\dot{z}_{j-\frac{1}{2}}\right)^2jdt-\frac{h}{2}\int_{0}^{T}\dot{z}_{N+\frac{1}{2}}\delta_x\dot{z}_{N+\frac{1}{2}}dt\\
		\nonumber 	&=-\frac{h}{2}\sum_{j=1}^{N}\int_{0}^{T}\left(\dot{z}_{j+\frac{1}{2}}\right)^2jdt+\frac{h}{2}\sum_{j=0}^{N-1}\int_{0}^{T}\left(\dot{z}_{j+\frac{1}{2}}\right)^2(j+1)dt-\frac{h}{2}\int_{0}^{T}\dot{z}_{N+\frac{1}{2}}\delta_x\dot{z}_{N+\frac{1}{2}}dt\\
\nonumber 			&=-\frac{h}{2}\sum_{j=0}^{N}\int_{0}^{T}\left(\dot{z}_{j+\frac{1}{2}}\right)^2jdt+\frac{h}{2}\sum_{j=0}^{N-1}\int_{0}^{T}\left(\dot{z}_{j+\frac{1}{2}}\right)^2jdt+\frac{h}{2}\sum_{j=0}^{N-1}\int_{0}^{T}\left(\dot{z}_{j+\frac{1}{2}}\right)^2dt-\frac{h}{2}\int_{0}^{T}\dot{z}_{N+\frac{1}{2}}\delta_x\dot{z}_{N+\frac{1}{2}}dt\\
	\nonumber 		&=-\frac{hN}{2}\int_{0}^{T}\left(\dot{z}_{N+\frac{1}{2}}\right)^2dt+\frac{h}{2}\sum_{j=0}^{N-1}\int_{0}^{T}\left(\dot{z}_{j+\frac{1}{2}}\right)^2dt-\frac{h}{2}\int_{0}^{T}\dot{z}_{N+\frac{1}{2}}\delta_x\dot{z}_{N+\frac{1}{2}}dt\\
	\nonumber 		&=\frac{(h-1)}{2}\int_{0}^{T}\left(\dot{z}_{N+\frac{1}{2}}\right)^2dt+\frac{h}{2}\sum_{j=0}^{N-1}\int_{0}^{T}\left(\dot{z}_{j+\frac{1}{2}}\right)^2dt-\frac{h}{2}\int_{0}^{T}\left(\frac{\dot{z}_{N+1}+\dot{z}_N}{2}\right)\left(\frac{\dot{z}_{N+1}-\dot{z}_N}{h}\right) dt\\
		\nonumber 	&=\frac{h}{2}\sum_{j=0}^{N}\int_{0}^{T}\left(\dot{z}_{j+\frac{1}{2}}\right)^2dt-\frac{1}{8}\int_{0}^{T}\left[3\left(\dot{z}_{N+1}\right)^2+2\dot{z}_{N+1}\dot{z}_N-\left(\dot{z}_N\right)^2\right]dt\\
		\nonumber 	&\geq \frac{h}{2}\sum_{j=0}^{N}\int_{0}^{T}\left(\dot{z}_{j+\frac{1}{2}}\right)^2dt-\frac{1}{8}\int_{0}^{T}\left[3\left(\dot{z}_{N+1}\right)^2+\left(\dot{z}_{N+1}\right)^2+\left(\dot{z}_{N}\right)^2-\left(\dot{z}_N\right)^2\right]dt\\
\label{denk6b}			&= \frac{h}{2}\sum_{j=0}^{N}\int_{0}^{T}\left(\dot{z}_{j+\frac{1}{2}}\right)^2dt-\frac{1}{2}\int_{0}^{T}\left(\dot{z}_{N+1}\right)^2dt.
		\end{align}
	Next, denote the first term of $I_2$ in \eqref{denk3} by $I_{2,1}$. Since $z_{j+1}-z_{j-1}=h(\delta_xz_{j+\frac{1}{2}}+\delta_xz_{j-\frac{1}{2}})$ and \\ $\delta_x^4 z_j=\frac{\delta_x^3z_{j+\frac{1}{2}}-\delta_x^3z_{j-\frac{1}{2}}}{h},$
		\begin{align}
			I_{2,1}&:=h\sum_{j=1}^{N}\int_{0}^{T}\delta_x^4z_jj\frac{z_{j+1}-z_{j-1}}{2}dt
			=\frac{h}{2}\sum_{j=1}^{N}\int_{0}^{T}j\left(\delta_x^3z_{j+\frac{1}{2}}-\delta_x^3z_{j-\frac{1}{2}}\right)(\delta_xz_{j+\frac{1}{2}}+\delta_xz_{j-\frac{1}{2}})dt\nonumber\\\nonumber
			&=\frac{h}{2}\sum_{j=1}^{N}\int_{0}^{T}j\delta_x^3z_{j+\frac{1}{2}}\delta_xz_{j+\frac{1}{2}}dt+\frac{h}{2}\sum_{j=1}^{N}\int_{0}^{T}j\delta_x^3z_{j+\frac{1}{2}}\delta_xz_{j-\frac{1}{2}}dt\\
			&\quad -\frac{h}{2}\sum_{j=1}^{N}\int_{0}^{T}j\delta_x^3z_{j-\frac{1}{2}}\delta_xz_{j+\frac{1}{2}}dt-\frac{h}{2}\sum_{j=1}^{N}\int_{0}^{T}j\delta_x^3z_{j-\frac{1}{2}}\delta_xz_{j-\frac{1}{2}}dt.\label{eq42}
		\end{align}
		Define $I_{2,2}$ by the first and the fourth terms  in \eqref{eq42} as the following
		\begin{align}
			I_{2,2}&:=\frac{h}{2}\sum_{j=1}^{N}\int_{0}^{T}j\delta_x^3z_{j+\frac{1}{2}}\delta_xz_{j+\frac{1}{2}}dt-\frac{h}{2}\sum_{j=1}^{N}\int_{0}^{T}j\delta_x^3z_{j-\frac{1}{2}}\delta_xz_{j-\frac{1}{2}}dt\\
			&=\frac{h}{2}\sum_{j=1}^{N}\int_{0}^{T}j\delta_x^3z_{j+\frac{1}{2}}\delta_xz_{j+\frac{1}{2}}dt-\frac{h}{2}\sum_{j=0}^{N-1}\int_{0}^{T}(j+1)\delta_x^3z_{j+\frac{1}{2}}\delta_xz_{j+\frac{1}{2}}dt\\
			&=\left(\frac{h}{2}\sum_{j=0}^{N-1}\int_{0}^{T}j\delta_x^3z_{j+\frac{1}{2}}\delta_xz_{j+\frac{1}{2}}dt+\frac{h}{2}\int_{0}^{T}N\delta_x^3z_{N+\frac{1}{2}}\delta_xz_{N+\frac{1}{2}}dt\right)\\
			&\quad-\frac{h}{2}\sum_{j=0}^{N-1}\int_{0}^{T}j\delta_x^3z_{j+\frac{1}{2}}\delta_xz_{j+\frac{1}{2}}dt-\frac{h}{2}\sum_{j=0}^{N-1}\int_{0}^{T}\delta_x^3z_{j+\frac{1}{2}}\delta_xz_{j+\frac{1}{2}}dt\\
			&=\frac{hN}{2}\int_{0}^{T}\delta_x^3z_{N+\frac{1}{2}}\delta_xz_{N+\frac{1}{2}}dt-\frac{h}{2}\sum_{j=0}^{N-1}\int_{0}^{T}\delta_x^3z_{j+\frac{1}{2}}\delta_xz_{j+\frac{1}{2}}dt\\
			&=\frac{(1-h)}{2}\int_{0}^{T}\delta_x^3z_{N+\frac{1}{2}}\delta_xz_{N+\frac{1}{2}}dt-\frac{h}{2}\sum_{j=0}^{N-1}\int_{0}^{T}\delta_x^3z_{j+\frac{1}{2}}\delta_xz_{j+\frac{1}{2}}dt\\
			\label{denk8}&=\frac{1}{2}\int_{0}^{T}\delta_x^3z_{N+\frac{1}{2}}\delta_xz_{N+\frac{1}{2}}dt-\frac{h}{2}\sum_{j=0}^{N}\int_{0}^{T}\delta_x^3z_{j+\frac{1}{2}}\delta_xz_{j+\frac{1}{2}}dt.
		\end{align}
		Substituting $I_{2,2}$ into $I_{2,1}$ yields
		\begin{align}
			\nonumber I_{2,1}&=\frac{1}{2}\int_{0}^{T}\delta_x^3z_{N+\frac{1}{2}}\delta_xz_{N+\frac{1}{2}}dt-\frac{h}{2}\sum_{j=0}^{N}\int_{0}^{T}\delta_x^3z_{j+\frac{1}{2}}\delta_xz_{j+\frac{1}{2}}dt\\
			\nonumber &\quad  +\frac{h}{2}\sum_{j=1}^{N}\int_{0}^{T}j\delta_x^3z_{j+\frac{1}{2}}\delta_xz_{j-\frac{1}{2}}dt-\frac{h}{2}\sum_{j=1}^{N}\int_{0}^{T}j\delta_x^3z_{j-\frac{1}{2}}\delta_xz_{j+\frac{1}{2}}dt\\
\nonumber &=\frac{1}{2}\int_{0}^{T}\delta_x^3z_{N+\frac{1}{2}}\delta_xz_{N+\frac{1}{2}}dt-\frac{h}{2}\sum_{j=0}^{N}\int_{0}^{T}\delta_x^3z_{j+\frac{1}{2}}\delta_xz_{j+\frac{1}{2}}dt\\ \nonumber
			&\quad+\frac{h}{2}\sum_{j=1}^{N}\int_{0}^{T}j\delta_x^3z_{j+\frac{1}{2}}\delta_xz_{j-\frac{1}{2}}dt-\frac{h}{2}\sum_{j=1}^{N}\int_{0}^{T}j\delta_x^3z_{j-\frac{1}{2}}\delta_xz_{j+\frac{1}{2}}dt\\
			\nonumber &=\frac{1}{2}\int_{0}^{T}\delta_x^3z_{N+\frac{1}{2}}\delta_xz_{N+\frac{1}{2}}dt-\frac{h}{2}\sum_{j=0}^{N}\int_{0}^{T}\delta_x^3z_{j+\frac{1}{2}}\delta_xz_{j+\frac{1}{2}}dt+I_{2,2}-\left(  \frac{h}{2}\sum_{j=1}^{N}\int_{0}^{T}j\delta_x^3z_{j+\frac{1}{2}}\left(\delta_xz_{j+\frac{1}{2}}-\delta_xz_{j-\frac{1}{2}}\right)dt\right)\\
\nonumber &\quad -\left( \frac{h}{2}\sum_{j=1}^{N}\int_{0}^{T}j\delta_x^3z_{j-\frac{1}{2}}\left(\delta_xz_{j+\frac{1}{2}}-\delta_xz_{j-\frac{1}{2}}\right)dt \right)\\
			\nonumber &=\int_{0}^{T}\delta_x^3z_{N+\frac{1}{2}}\delta_xz_{N+\frac{1}{2}}dt-h\sum_{j=0}^{N}\int_{0}^{T}\delta_x^3z_{j+\frac{1}{2}}\delta_xz_{j+\frac{1}{2}}dt -\left(  \frac{h}{2}\sum_{j=1}^{N}\int_{0}^{T}hj\delta_x^3z_{j+\frac{1}{2}}\delta_x^2z_jdt\right)\\
	\label{denk9}	&\quad-\left( \frac{h}{2}\sum_{j=1}^{N}\int_{0}^{T}hj\delta_x^3z_{j-\frac{1}{2}}\delta_x^2z_jdt \right).
		\end{align}
	Now, substitute $\delta_x^3z_{j+\frac{1}{2}}=\frac{\delta_x^2z_{j+1}-\delta_x^2z_{j}}{h}$ in \eqref{denk9}, and use  $\delta_xz_{-\frac{1}{2}}=\delta_x^2z_{N+1}=0$
		\begin{align}
	\nonumber 		I_{2,1}=&\int_{0}^{T}\delta_x^3z_{N+\frac{1}{2}}\delta_xz_{N+\frac{1}{2}}dt-h\sum_{j=0}^{N}\int_{0}^{T}\frac{\delta_x^2z_{j+1}-\delta_x^2z_{j}}{h}\delta_xz_{j+\frac{1}{2}}dt\\
	\nonumber 		&-\left(  \frac{h}{2}\sum_{j=1}^{N}\int_{0}^{T}hj\frac{\delta_x^2z_{j+1}-\delta_x^2z_{j}}{h}\delta_x^2z_jdt\right)-\left( \frac{h}{2}\sum_{j=1}^{N}\int_{0}^{T}hj\frac{\delta_x^2z_{j}-\delta_x^2z_{j-1}}{h}\delta_x^2z_jdt \right)\\
\nonumber 			=&\int_{0}^{T}\delta_x^3z_{N+\frac{1}{2}}\delta_xz_{N+\frac{1}{2}}dt-\sum_{j=1}^{N+1}\int_{0}^{T}\delta_x^2z_{j}\delta_xz_{j-\frac{1}{2}}dt+\sum_{j=0}^{N}\int_{0}^{T}\delta_x^2z_{j}\delta_xz_{j+\frac{1}{2}}dt\\
	\nonumber 		&-  \frac{h}{2}\sum_{j=1}^{N}\int_{0}^{T}j\delta_x^2z_{j+1}\delta_x^2z_jdt+ \frac{h}{2}\sum_{j=0}^{N-1}\int_{0}^{T}(j+1)\delta_x^2z_{j}\delta_x^2z_{j+1}dt\\
\nonumber 			=&\int_{0}^{T}\delta_x^3z_{N+\frac{1}{2}}\delta_xz_{N+\frac{1}{2}}dt-\sum_{j=0}^{N}\int_{0}^{T}\delta_x^2z_{j}\delta_xz_{j-\frac{1}{2}}dt+\underbrace{\int_{0}^{T}\delta_x^2z_{0}\delta_xz_{-\frac{1}{2}}dt  -\int_{0}^{T}\delta_x^2z_{N+1}\delta_xz_{N+\frac{1}{2}}dt}_{=0}\\
	\nonumber 		&+\sum_{j=0}^{N}\int_{0}^{T}\delta_x^2z_{j}\delta_xz_{j+\frac{1}{2}}dt-\underbrace{\frac{h(N+1)}{2}\int_{0}^{T}\delta_x^2z_{N}\delta_x^2z_{N+1}dt}_{=0}\\
\nonumber 			&-  \frac{h}{2}\sum_{j=0}^{N}\int_{0}^{T}j\delta_x^2z_{j+1}\delta_x^2z_jdt+ \frac{h}{2}\sum_{j=0}^{N}\int_{0}^{T}j\delta_x^2z_{j}\delta_x^2z_{j+1}dt+\frac{h}{2}\sum_{j=0}^{N}\int_{0}^{T}\delta_x^2z_{j}\delta_x^2z_{j+1}dt\\
	\label{denk10}		=&\int_{0}^{T}\delta_x^3z_{N+\frac{1}{2}}\delta_xz_{N+\frac{1}{2}}dt+h\sum_{j=0}^{N}\int_{0}^{T}\left(\delta_x^2z_{j}\right)^2dt+\frac{h}{2}\sum_{j=0}^{N}\int_{0}^{T}\delta_x^2z_{j}\delta_x^2z_{j+1}dt.
		\end{align}

		Now, substitute \eqref{denk10} into $I_2$ in \eqref{denk3} to obtain
		\begin{equation}\label{i2}
			I_2=h\sum_{j=0}^{N}\int_{0}^{T}\left(\delta_x^2z_{j}\right)^2dt+\frac{h}{2}\sum_{j=0}^{N}\int_{0}^{T}\delta_x^2z_{j}\delta_x^2z_{j+1}dt.
		\end{equation}
		Note also that  by \eqref{denk6b}
		\begin{equation}\label{i1}
			I_1\geq\frac{h}{2}\sum_{j=0}^{N}\int_{0}^{T}\left(\dot{z}_{j+\frac{1}{2}}\right)^2dt-\frac{1}{2}\int_{0}^{T}\left(\dot{z}_{N+1}\right)^2dt.
		\end{equation}
		Thus, by the inequality $2ab\geq -a^2-b^2$,  \eqref{i2} and \eqref{i1}, the sum of the first three terms $L_{h,1}(t)|_0^T+I_1+I_2$ can be estimated by the  energy as the following
		\begin{align}
			L_{h,1}(t)|_0^T+I_1+I_2&\geq -2E_h(0)+\frac{h}{2}\sum_{j=0}^{N}\int_{0}^{T}\left(\dot{z}_{j+\frac{1}{2}}\right)^2dt-\frac{1}{2}\int_{0}^{T}\left(\dot{z}_{N+1}\right)^2dt\nonumber\\
			&+h\sum_{j=0}^{N}\int_{0}^{T}\left(\delta_x^2z_{j}\right)^2dt+\frac{h}{2}\sum_{j=0}^{N}\int_{0}^{T}\delta_x^2z_{j}\delta_x^2z_{j+1}dt\nonumber\\
			\geq& -2E_h(0)+\frac{h}{2}\sum_{j=0}^{N}\int_{0}^{T}\left(\dot{z}_{j+\frac{1}{2}}\right)^2dt-\frac{1}{2}\int_{0}^{T}\left(\dot{z}_{N+1}\right)^2dt\nonumber\\
			&+\frac{3h}{4}\sum_{j=0}^{N}\int_{0}^{T}\left(\delta_x^2z_{j}\right)^2dt-\frac{h}{4}\sum_{j=1}^{N+1}\int_{0}^{T}\left(\delta_x^2z_{j}\right)^2dt\nonumber\\
			=& -2E_h(0)-\frac{1}{2}\int_{0}^{T}\left(\dot{z}_{N+1}\right)^2dt+\frac{h}{2}\sum_{j=0}^{N}\int_{0}^{T}\left(\dot{z}_{j+\frac{1}{2}}\right)^2dt+\frac{h}{2}\sum_{j=0}^{N}\int_{0}^{T}\left(\delta_x^2z_{j}\right)^2dt.\label{firsttermssum}
		\end{align}
		Next, consider the first term of $I_3$ in \eqref{denk3} and denote it by
		\begin{align*}
			I_{3,1}:=&-\frac{B}{4}\sum_{j=1}^{N}\int_{0}^{T}j(\phi_{j+1}-\phi_{j-1})y_jdt\\
			=&-\frac{B}{4}\sum_{j=1}^{N}\int_{0}^{T}j\phi_{j+1}y_jdt+\frac{B}{4}\sum_{j=0}^{N-1}\int_{0}^{T}(j+1)\phi_{j}y_{j+1}dt\\
			=&-\frac{B}{4}\sum_{j=1}^{N}\int_{0}^{T}j\phi_{j+1}y_jdt+\frac{B}{4}\sum_{j=1}^{N}\int_{0}^{T}j\phi_{j}y_{j+1}dt+\frac{B}{4}\sum_{j=1}^{N}\int_{0}^{T}\phi_{j}y_{j+1}dt\\
			&+\underbrace{\frac{B}{4}\int_{0}^{T}\phi_{0}y_{1}dt}_{=0}-\frac{B(N+1)}{4}\int_{0}^{T}\phi_{N}y_{N+1}dt.
		\end{align*}
		Adding and subtracting appropriate terms, by the boundary conditions $\phi_N=\phi_{N+1}$, $\phi_0=0$, and $y_{N+1}=y_N,$
		\begin{align*}
			I_{3,1}=&-\frac{B}{4}\sum_{j=1}^{N}\int_{0}^{T}\phi_{j+1}jy_jdt+\frac{B}{4}\sum_{j=1}^{N}\int_{0}^{T}\phi_{j}jy_{j+1}dt+\frac{B}{4}\sum_{j=1}^{N}\int_{0}^{T}\phi_{j}y_{j+1}dt-\frac{B(N+1)}{4}\int_{0}^{T}\phi_{N}y_{N+1}dt\\
			&+\underbrace{\frac{B}{4}\sum_{j=1}^{N}\int_{0}^{T}\phi_{j+1}jy_{j+1}dt-\frac{B}{4}\sum_{j=1}^{N}\int_{0}^{T}\phi_{j+1}jy_{j+1}dt}_{=0}+\underbrace{\frac{B}{4}\sum_{j=1}^{N}\int_{0}^{T}\phi_{j}jy_{j}dt-\frac{B}{4}\sum_{j=1}^{N}\int_{0}^{T}\phi_{j}jy_{j}dt}_{=0}\\			=&\frac{B}{4}\sum_{j=1}^{N}\int_{0}^{T}j\left(\phi_{j+1}y_{j+1}-\phi_{j+1}y_{j}+\phi_{j}y_{j+1}-\phi_{j}y_{j} \right)dt-\frac{B(N+1)}{4}\int_{0}^{T}\phi_{N}y_{N+1}dt\\
			&+\frac{B}{4}\sum_{j=1}^{N}\int_{0}^{T}\phi_{j}jy_{j}dt-\frac{B}{4}\sum_{j=1}^{N}\int_{0}^{T}\phi_{j+1}jy_{j+1}dt+\frac{B}{4}\sum_{j=1}^{N}\int_{0}^{T}\phi_{j}y_{j+1}dt\\
			=&\frac{Bh}{2}\sum_{j=1}^{N}\int_{0}^{T}j\phi_{j+\frac{1}{2}}\delta_xy_{j+\frac{1}{2}}dt +\frac{B}{4}\sum_{j=1}^{N}\int_{0}^{T}\phi_{j}(y_{j+1}+y_j)dt-\frac{B(2N+1)}{4}\int_{0}^{T}\phi_{N}y_{N}dt.
		\end{align*}
		\normalsize
		Consider the first term of $I_{3,1}$ together with \eqref{kdef} and the boundary conditions $k_0=\delta_xk_{N+\frac{1}{2}}=0$ so that
		\begin{align*}
			I_{3,2}:=&\frac{Bh}{2}\sum_{j=1}^{N}\int_{0}^{T}j\phi_{j+\frac{1}{2}}\delta_xy_{j+\frac{1}{2}}dt\\
			=\frac{B}{4}&\sum_{j=1}^{N}\int_{0}^{T}j\frac{B}{2Ph}\left(\left(\bm{A}_h\vec{k}\right)_{j+1}+\left(\bm{A}_h\vec{k}\right)_{j}\right) \left(\left (\frac{C}{P}-\frac{h^2}{4}\right)\left(\left(\bm{A}_h\vec{k}\right)_{j+1}-\left(\bm{A}_h\vec{k}\right)_{j}\right)+(k_{j+1}-k_j)\right)dt\\
			&+\frac{B^2}{8Ph}\sum_{j=1}^{N}\int_{0}^{T}j\left(\left(\bm{A}_h\vec{k}\right)_{j+1}+\left(\bm{A}_h\vec{k}\right)_{j}\right)(k_{j+1}-k_j)dt\\
			=&\frac{B^2}{8Ph}\left(\frac{C}{P}-\frac{h^2}{4}\right)\sum_{j=2}^{N+1}\int_{0}^{T}(j-1)\left(\bm{A}_h\vec{k}\right)_{j}^2dt-\frac{B^2}{8Ph}\left(\frac{C}{P}-\frac{h^2}{4}\right)\sum_{j=1}^{N}\int_{0}^{T}j\left(\bm{A}_h\vec{k}\right)_{j}^2dt\\
			&+\frac{B^2}{8Ph}\sum_{j=1}^{N}\int_{0}^{T}j\left(\left(\bm{A}_h\vec{k}\right)_{j+1}+\left(\bm{A}_h\vec{k}\right)_{j}\right)(k_{j+1}-k_j)dt\\
			=&-\frac{B^2}{8Ph}\left(\frac{C}{P}-\frac{h^2}{4}\right)\sum_{j=1}^{N}\int_{0}^{T}\left(\bm{A}_h\vec{k}\right)_{j}^2dt+\frac{B^2N}{8Ph}\left(\frac{C}{P}-\frac{h^2}{4}\right)\int_{0}^{T}\left(\bm{A}_h\vec{k}\right)_{N+1}^2dt\\
			&+\frac{B^2}{8Ph}\sum_{j=1}^{N}\int_{0}^{T}j\left(\left(\bm{A}_h\vec{k}\right)_{j+1}+\left(\bm{A}_h\vec{k}\right)_{j}\right)(k_{j+1}-k_j)dt
		\end{align*}
		\normalsize
		where
		\begin{align}
			\frac{B^2}{8Ph}&\sum_{j=1}^{N}\int_{0}^{T}j\left(\bm{A}_h\vec{k}\right)_{j}(k_{j+1}-k_j)dt=-\frac{B^2}{8Ph}\sum_{j=1}^{N}\int_{0}^{T}j\delta_x^2k_{j}(k_{j+1}-k_j)dt\nonumber\\
			=&-\frac{B^2}{8Ph}\sum_{j=1}^{N}\int_{0}^{T}j(\delta_xk_{j+\frac{1}{2}}-\delta_xk_{j-\frac{1}{2}})\delta_xk_{j+\frac{1}{2}dt}\nonumber\\
			=&-\frac{B^2}{8Ph}\sum_{j=0}^{N-1}\int_{0}^{T}j\delta_xk_{j+\frac{1}{2}}\delta_xk_{j+\frac{1}{2}}dt+\frac{B^2}{8Ph}\sum_{j=0}^{N-1}\int_{0}^{T}(j+1)\delta_xk_{j+\frac{1}{2}}\delta_xk_{j+\frac{3}{2}}dt\\
&\quad\quad -\underbrace{\frac{B^2N}{8Ph}\int_{0}^{T}\left(\delta_xk_{N+\frac{1}{2}}\right)^2dt}_{=0 \text{ by the B.C.'s}}\nonumber\\
			=&-\frac{B^2}{8Ph}\sum_{j=1}^{N}\int_{0}^{T}j\left(\bm{A}_h\vec{k}\right)_{j+1}(k_{j+1}-k_j)dt+\frac{B^2}{8Ph}\sum_{j=1}^{N}\int_{0}^{T}\delta_xk_{j-\frac{1}{2}}\delta_xk_{j+\frac{1}{2}}dt\label{ibp}
		\end{align}
	by $k_{N+1}=k_N.$	By collecting similar terms to the same side of \eqref{ibp}
		\begin{eqnarray}
			 \label{secondpart} &&	\frac{B^2}{8Ph}\sum_{j=1}^{N}\int_{0}^{T}j \left(\left(\bm{A}_h\vec{k}\right)_{j+1}+\left(\bm{A}_h\vec{k}\right)_{j}\right) (k_{j+1}-k_j)dt =\frac{B^2}{8Ph}\sum_{j=1}^{N}\int_{0}^{T}\delta_xk_{j-\frac{1}{2}}\delta_xk_{j+\frac{1}{2}}dt,
					\end{eqnarray}
		and substituting \eqref{secondpart} into $I_{3,2}$ above
		\begin{equation}\label{i31}
			\begin{split}
				I_{3,2}=&-\frac{B^2}{8Ph}\left(\frac{C}{P}-\frac{h^2}{4}\right)\sum_{j=1}^{N}\int_{0}^{T}\left(\bm{A}_h\vec{k}\right)_{j}^2dt+\frac{B^2}{8Ph}\sum_{j=1}^{N}\int_{0}^{T}\delta_xk_{j-\frac{1}{2}}\delta_xk_{j+\frac{1}{2}}dt\\
				&+\frac{B^2N}{8Ph}\left(\frac{C}{P}-\frac{h^2}{4}\right)\int_{0}^{T}\left(\bm{A}_h\vec{k}\right)_{N+1}^2dt.
			\end{split}.
		\end{equation}
		Now, consider the second term of $I_{3,1}$ and use definitions in \eqref{kdef}
		\begin{align*}
			I_{3,3}:=&\frac{B}{4}\sum_{j=1}^{N}\int_{0}^{T}\phi_{j}(y_{j+1}+y_j)dt\\
			=&\frac{B}{4}\sum_{j=1}^{N}\int_{0}^{T}\frac{B}{2Ph}\left(\bm{A}_h\vec{k}\right)_{j}\left( \left(\frac{C}{P}-\frac{h^2}{4}\right)\left(\left(\bm{A}_h\vec{k}\right)_{j+1}+\left(\bm{A}_h\vec{k}\right)_{j}\right)+(k_{j+1}+k_j)\right)dt\\
			=&\frac{B^2}{8Ph}\left(\frac{C}{P}-\frac{h^2}{4}\right)\sum_{j=1}^{N}\int_{0}^{T}\left(\bm{A}_h\vec{k}\right)_{j}^2dt+\frac{B^2}{8Ph}\left(\frac{C}{P}-\frac{h^2}{4}\right)\sum_{j=1}^{N}\int_{0}^{T}\left(\bm{A}_h\vec{k}\right)_{j}\left(\bm{A}_h\vec{k}\right)_{j+1}dt\\
			&+\frac{B^2}{8Ph}\sum_{j=1}^{N}\int_{0}^{T}\left(\bm{A}_h\vec{k}\right)_{j}k_{j+1}dt+\frac{B^2}{8Ph}\sum_{j=1}^{N}\int_{0}^{T}\left(\bm{A}_h\vec{k}\right)_{j}k_{j}dt\\		=&\frac{B^2}{8Ph}\left(\frac{C}{P}-\frac{h^2}{4}\right)\sum_{j=1}^{N}\int_{0}^{T}\left(\bm{A}_h\vec{k}\right)_{j}^2+\frac{B^2}{8Ph}\left(\frac{C}{P}-\frac{h^2}{4}\right)\sum_{j=1}^{N}\int_{0}^{T}\left(\bm{A}_h\vec{k}\right)_{j}\left(\bm{A}_h\vec{k}\right)_{j+1}dt\\
			&-\frac{B^2}{8Ph}\sum_{j=1}^{N}\int_{0}^{T}\delta_x^2k_jk_{j+1}dt-\frac{B^2}{8Ph}\sum_{j=1}^{N}\int_{0}^{T}\delta_x^2k_jk_{j}dt\\
			=&\frac{B^2}{8Ph}\left(\frac{C}{P}-\frac{h^2}{4}\right)\sum_{j=1}^{N}\int_{0}^{T}\left(\bm{A}_h\vec{k}\right)_{j}^2dt+\frac{B^2}{8Ph}\left(\frac{C}{P}-\frac{h^2}{4}\right)\sum_{j=1}^{N}\int_{0}^{T}\left(\bm{A}_h\vec{k}\right)_{j}\left(\bm{A}_h\vec{k}\right)_{j+1}dt\\
			&-\frac{B^2}{8Ph^2}\sum_{j=1}^{N}\int_{0}^{T}\left(\delta_xk_{j+\frac{1}{2}}-\delta_xk_{j-\frac{1}{2}}\right) k_{j+1}dt-\frac{B^2}{8Ph^2}\sum_{j=1}^{N}\int_{0}^{T}\left(\delta_xk_{j+\frac{1}{2}}-\delta_xk_{j-\frac{1}{2}}\right)k_{j}dt\\
			=&\frac{B^2}{8Ph}\left(\frac{C}{P}-\frac{h^2}{4}\right)\sum_{j=1}^{N}\int_{0}^{T}\left(\bm{A}_h\vec{k}\right)_{j}^2dt+\frac{B^2}{8Ph}\left(\frac{C}{P}-\frac{h^2}{4}\right)\sum_{j=1}^{N}\int_{0}^{T}\left(\bm{A}_h\vec{k}\right)_{j}\left(\bm{A}_h\vec{k}\right)_{j+1dt}\\
			&-\frac{B^2}{8Ph^2}\sum_{j=2}^{N+1}\int_{0}^{T}\delta_xk_{j-\frac{1}{2}}k_{j}dt+\frac{B^2}{8Ph^2}\sum_{j=1}^{N}\int_{0}^{T}\delta_xk_{j-\frac{1}{2}}k_{j+1}dt\\
			&-\frac{B^2}{8Ph^2}\sum_{j=1}^{N}\int_{0}^{T}\delta_xk_{j+\frac{1}{2}}k_{j}dt+\frac{B^2}{8Ph^2}\sum_{j=0}^{N-1}\int_{0}^{T}\delta_xk_{j+\frac{1}{2}}k_{j+1}dt\\
			=&\frac{B^2}{8Ph}\left(\frac{C}{P}-\frac{h^2}{4}\right)\sum_{j=1}^{N}\int_{0}^{T}\left(\bm{A}_h\vec{k}\right)_{j}^2dt+\frac{B^2}{8Ph}\left(\frac{C}{P}-\frac{h^2}{4}\right)\sum_{j=1}^{N}\int_{0}^{T}\left(\bm{A}_h\vec{k}\right)_{j}\left(\bm{A}_h\vec{k}\right)_{j+1}dt\\
			&+\frac{B^2}{8Ph}\sum_{j=1}^{N}\int_{0}^{T}\delta_xk_{j-\frac{1}{2}}\delta_xk_{j+\frac{1}{2}}dt-\underbrace{\frac{B^2}{8Ph^2}\int_{0}^{T}\delta_xk_{N+\frac{1}{2}}k_{N+1}dt}_{=0}+\frac{B^2}{8Ph^2}\int_{0}^{T}\delta_xk_{\frac{1}{2}}k_{1}dt\\
			&+\frac{B^2}{8Ph}\sum_{j=1}^{N}\int_{0}^{T}\delta_xk_{j+\frac{1}{2}}\delta_xk_{j+\frac{1}{2}}dt-\underbrace{\frac{B^2}{8Ph^2}\int_{0}^{T}\delta_xk_{N+\frac{1}{2}}k_{N+1}dt}_{=0}+\frac{B^2}{8Ph^2}\int_{0}^{T}\delta_xk_{\frac{1}{2}}k_{1}dt\\
			=&\frac{B^2}{8Ph}\left(\frac{C}{P}-\frac{h^2}{4}\right)\sum_{j=1}^{N}\int_{0}^{T}\left(\bm{A}_h\vec{k}\right)_{j}^2dt+\frac{B^2}{8Ph}\left(\frac{C}{P}-\frac{h^2}{4}\right)\sum_{j=1}^{N}\int_{0}^{T}\left(\bm{A}_h\vec{k}\right)_{j}\left(\bm{A}_h\vec{k}\right)_{j+1}dt\\
			&+\frac{B^2}{8Ph}\sum_{j=1}^{N}\int_{0}^{T}\delta_xk_{j-\frac{1}{2}}\delta_xk_{j+\frac{1}{2}}dt+\frac{B^2}{8Ph}\sum_{j=1}^{N}\int_{0}^{T}\delta_xk_{j+\frac{1}{2}}\delta_xk_{j+\frac{1}{2}}dt+\frac{B^2}{4Ph^3}\int_{0}^{T}k_{1}^2dt.		
		\end{align*}
		Recalling
		\begin{equation}\label{i3sum} I_3=I_{3,1}+B\int_{0}^{T}\phi_N\delta_xz_{N+\frac{1}{2}}dt=I_{3,2}+I_{3,3}-\frac{B(2N+1)}{4}\int_{0}^{T}\phi_{N}y_{N}dt+B\int_{0}^{T}\phi_N\delta_xz_{N+\frac{1}{2}}dt,
		\end{equation}
		\normalsize
		and substituting $I_{3,2}$ and $I_{3,3}$ into \eqref{i3sum} lead to
		\begin{equation}\label{i3exp}
			\begin{split}
				I_3=&\frac{B^2}{8Ph}\left(\frac{C}{P}-\frac{h^2}{4}\right)\sum_{j=1}^{N}\int_{0}^{T}\left(\bm{A}_h\vec{k}\right)_{j}\left(\bm{A}_h\vec{k}\right)_{j+1}dt+\frac{B^2}{4Ph}\sum_{j=1}^{N}\int_{0}^{T}\delta_xk_{j-\frac{1}{2}}\delta_xk_{j+\frac{1}{2}}dt+\frac{B}{4}\int_{0}^{T}\phi_Ny_{N}dt\\
				&+\frac{B^2N}{8Ph}\left(\frac{C}{P}-\frac{h^2}{4}\right)\int_{0}^{T}\left(\bm{A}_h\vec{k}\right)_{N+1}^2dt+\frac{B^2}{8Ph}\sum_{j=1}^{N}\int_{0}^{T}\delta_xk_{j+\frac{1}{2}}\delta_xk_{j+\frac{1}{2}}dt+\frac{B^2}{4Ph^3}\int_{0}^{T}k_{1}^2dt
			\end{split}
		\end{equation}	
		where boundary terms are simplified by $z_N=\frac{z_{N+1}+z_{N-1}}{2}$ as the following
		\begin{align*}
			B\int_{0}^{T}\phi_N\delta_xz_{N+\frac{1}{2}}dt-\frac{B(2N+1)}{4}\int_{0}^{T}\phi_{N}y_{N}dt&=\frac{B}{h}\int_{0}^{T}\phi_N\left(z_{N+1}-z_N\right)dt-\frac{B(2N+1)}{4}\int_{0}^{T}\phi_{N}y_{N}dt\\
			&=\frac{B(N+1)}{2}\int_{0}^{T}\phi_Ny_Ndt-\frac{B(2N+1)}{4}\int_{0}^{T}\phi_{N}y_{Ndt}\\
			&=\frac{B}{4}\int_{0}^{T}\phi_Ny_{N}dt.
		\end{align*}
		Next, consider following terms in \eqref{i3exp}
		\begin{align}
			&\frac{B^2}{4Ph}\sum_{j=1}^{N}\int_{0}^{T}\delta_xk_{j-\frac{1}{2}}\delta_xk_{j+\frac{1}{2}}dt+\frac{B^2}{8Ph}\sum_{j=1}^{N}\int_{0}^{T}\delta_xk_{j+\frac{1}{2}}\delta_xk_{j+\frac{1}{2}}dt\nonumber\\
			&=\frac{B^2}{4Ph}\sum_{j=1}^{N}\int_{0}^{T}\delta_xk_{j-\frac{1}{2}}\delta_xk_{j+\frac{1}{2}}dt+\frac{B^2}{16Ph}\sum_{j=1}^{N}\int_{0}^{T}\left(\delta_xk_{j+\frac{1}{2}}\right)^2dt+\frac{B^2}{16Ph}\sum_{j=2}^{N+1}\int_{0}^{T}\left(\delta_xk_{j-\frac{1}{2}}\right)^2dt\nonumber\\
			&=\frac{B^2}{8Ph}\sum_{j=1}^{N}\int_{0}^{T}\delta_xk_{j-\frac{1}{2}}\delta_xk_{j+\frac{1}{2}}dt+\frac{B^2}{16Ph}\sum_{j=2}^{N}\int_{0}^{T}\left(\delta_xk_{j+\frac{1}{2}}+\delta_xk_{j-\frac{1}{2}}\right)^2dt\nonumber\\
			&\quad +\frac{B^2}{16Ph}\int_{0}^{T}\left(\delta_xk_{N+\frac{1}{2}}\right)^2dt+\frac{B^2}{16Ph}\int_{0}^{T}\left(\delta_xk_{\frac{3}{2}}\right)^2dt+\frac{B^2}{8Ph}\int_{0}^{T}\delta_xk_{\frac{1}{2}}\delta_xk_{\frac{3}{2}}dt\nonumber\\
			&=\frac{B^2}{8Ph}\sum_{j=1}^{N}\int_{0}^{T}\delta_xk_{j-\frac{1}{2}}\delta_xk_{j+\frac{1}{2}}dt+\frac{B^2}{16Ph}\sum_{j=1}^{N}\int_{0}^{T}\left(\delta_xk_{j+\frac{1}{2}}+\delta_xk_{j-\frac{1}{2}}\right)^2dt-\frac{B^2}{16Ph}\int_{0}^{T}\left(\delta_xk_{\frac{1}{2}}\right)^2dt.\nonumber
		\end{align}
		and substitute them back into $I_3$
		\begin{equation*}
			\begin{split}
				I_3=&\frac{B^2}{8Ph}\left(\frac{C}{P}-\frac{h^2}{4}\right)\sum_{j=1}^{N}\int_{0}^{T}\left(\bm{A}_h\vec{k}\right)_{j}\left(\bm{A}_h\vec{k}\right)_{j+1}dt+\frac{B^2}{8Ph}\sum_{j=1}^{N}\int_{0}^{T}\delta_xk_{j-\frac{1}{2}}\delta_xk_{j+\frac{1}{2}}dt+\frac{B}{4}\int_{0}^{T}\phi_Ny_{N}dt\\
				&+\frac{B^2}{16Ph}\sum_{j=1}^{N}\int_{0}^{T}\left(\delta_xk_{j+\frac{1}{2}}dt+\delta_xk_{j-\frac{1}{2}}\right)^2dt+\frac{B^2N}{8Ph}\left(\frac{C}{P}-\frac{h^2}{4}\right)\int_{0}^{T}\left(\bm{A}_h\vec{k}\right)_{N+1}^2dt+\frac{3B^2}{16Ph^3}\int_{0}^{T}k_{1}^2dt.	
			\end{split}
		\end{equation*}
		where the boundary condition $k_0=0$ is used for $\delta_xk_{\frac{1}{2}}=\frac{k_1}{h}$.
Finally, we have an estimate for $I_3$
		\begin{equation}\label{i3final}
			\begin{split}
				I_3\geq&\frac{B^2}{8Ph}\left(\frac{C}{P}-\frac{h^2}{4}\right)\sum_{j=1}^{N}\int_{0}^{T}\left(\bm{A}_h\vec{k}\right)_{j}\left(\bm{A}_h\vec{k}\right)_{j+1}dt+\frac{B^2}{8Ph}\sum_{j=1}^{N}\int_{0}^{T}\delta_xk_{j-\frac{1}{2}}\delta_xk_{j+\frac{1}{2}}dt+\frac{B}{4}\int_{0}^{T}\phi_Ny_{N}dt.	
			\end{split}
		\end{equation}
		Now, it is time for the second multiplier.  Multiply the first equation in \eqref{clampeddiscrete} by $h\frac{z_j}{4},$ the third (boundary) equation in \eqref{clampeddiscrete} by $\frac{z_{N+1}}{4},$ respectively, and integrate by parts with respect to $t$ over $[0,T]$ to obtain
		\begin{align*}
			0=&\frac{h}{8}\sum_{j=1}^{N}\int_{0}^{T}\left(\ddot{z}_{j+\frac{1}{2}}+\ddot{z}_{j-\frac{1}{2}}\right)z_{j}dt+\frac{h}{4}\sum_{j=1}^{N}\int_{0}^{T}\delta_x^4z_jz_{j}dt-\frac{B}{8}\sum_{j=1}^{N}\int_{0}^{T}(\phi_{j+1}-\phi_{j-1})z_{j}dt\\
			=&\frac{h}{8}\sum_{j=1}^{N}\left.\left(\dot{z}_{j+\frac{1}{2}}+\dot{z}_{j-\frac{1}{2}}\right)z_{j+1}\right|_0^T-\frac{h}{8}\sum_{j=1}^{N}\int_{0}^{T}\left(\dot{z}_{j+\frac{1}{2}}+\dot{z}_{j-\frac{1}{2}}\right)\dot{z}_{j+1}dt\\
			&+\frac{h}{4}\sum_{j=1}^{N}\int_{0}^{T}\delta_x^4z_jz_{j}dt-\frac{B}{8}\sum_{j=1}^{N}\int_{0}^{T}(\phi_{j+1}-\phi_{j-1})z_{j}dt,
		\end{align*}
		\begin{align*}
			0=&\frac{h}{8}\int_{0}^{T}\ddot{z}_{N+\frac{1}{2}}z_{N+1}dt-\frac{1}{4}\int_{0}^{T}\delta_x^3z_{N+\frac{1}{2}}z_{N+1}dt+\frac{B}{4}\int_{0}^{T}\phi_Nz_{N+1}dt\\
			=&\left.\frac{h}{8}\dot{z}_{N+\frac{1}{2}}z_{N+1}\right|_0^T-\frac{h}{8}\int_{0}^{T}\dot{z}_{N+\frac{1}{2}}\dot{z}_{N+1}dt-\frac{1}{4}\int_{0}^{T}\delta_x^3z_{N+\frac{1}{2}}z_{N+1}dt+\frac{B}{4}\int_{0}^{T}\phi_Nz_{N+1}dt.
		\end{align*}
		Define
		\begin{align*}
			L_{h,2}(t)&:=\frac{h}{8}\sum_{j=1}^{N}\left(\dot{z}_{j+\frac{1}{2}}+\dot{z}_{j-\frac{1}{2}}\right)z_{j}+\frac{h}{8}\dot{z}_{N+\frac{1}{2}}z_{N+1},\\
			I_4&:=-\frac{h}{8}\sum_{j=1}^{N}\int_{0}^{T}\left(\dot{z}_{j+\frac{1}{2}}+\dot{z}_{j-\frac{1}{2}}\right)\dot{z}_{j}dt-\frac{h}{8}\int_{0}^{T}\dot{z}_{N+\frac{1}{2}}\dot{z}_{N+1}dt,\\
			I_5&:=\frac{h}{4}\sum_{j=1}^{N}\int_{0}^{T}\delta_x^4z_jz_{j}dt-\frac{1}{4}\int_{0}^{T}\delta_x^3z_{N+\frac{1}{2}}z_{N+1}dt,\\
			I_6&:=-\frac{B}{8}\sum_{j=1}^{N}\int_{0}^{T}(\phi_{j+1}-\phi_{j-1})z_{j}dt+\frac{B}{4}\int_{0}^{T}\phi_Nz_{N+1}dt.
		\end{align*}
		such that $L_{h,2}(t)|_0^T+I_4+I_5+I_6=0.$
		First of all, the Young's inequality with  \cref{poincare} imply
		\begin{align*}
			L_{h,2}(t)&=\frac{h}{8}\sum_{j=1}^{N}\left(\dot{z}_{j+\frac{1}{2}}+\dot{z}_{j-\frac{1}{2}}\right)z_{j}+\frac{h}{8}\dot{z}_{N+\frac{1}{2}}z_{N+1}\\
			&\leq \frac{h}{16}\sum_{j=1}^{N}\left(\dot{z}_{j+\frac{1}{2}}+\dot{z}_{j-\frac{1}{2}}\right)^2+\frac{h}{16}\sum_{j=1}^{N}z_{j}^2+\frac{h}{16}\left(\dot{z}_{N+\frac{1}{2}}\right)^2+\frac{h}{16}\left(z_{N+1}\right)^2\\
			&\leq  \frac{h}{8}\sum_{j=1}^{N}\dot{z}_{j+\frac{1}{2}}^2+\frac{h}{8}\sum_{j=1}^{N}\dot{z}_{j-\frac{1}{2}}^2+\frac{h}{16}\sum_{j=1}^{N}z_{j}^2+\frac{h}{16}\left(\dot{z}_{N+\frac{1}{2}}\right)^2+\frac{h}{16}\left(z_{N+1}\right)^2\\
			&=  \frac{h}{8}\sum_{j=1}^{N}\dot{z}_{j+\frac{1}{2}}^2+\frac{h}{8}\sum_{j=0}^{N-1}\dot{z}_{j+\frac{1}{2}}^2+\frac{h}{16}\sum_{j=1}^{N}z_{j}^2+\frac{h}{16}\left(\dot{z}_{N+\frac{1}{2}}\right)^2+\frac{h}{16}\left(z_{N+1}\right)^2\\
			&\leq  \frac{h}{8}\sum_{j=1}^{N}\dot{z}_{j+\frac{1}{2}}^2+\frac{h}{8}\sum_{j=0}^{N}\dot{z}_{j+\frac{1}{2}}^2+\frac{h}{16}\sum_{j=1}^{N}z_{j}^2+\frac{h}{16}\left(z_{N+1}\right)^2\\
			&\leq \frac{h}{4}\sum_{j=0}^{N}\dot{z}_{j+\frac{1}{2}}^2+\frac{h}{16}\sum_{j=1}^{N}\left(\delta_x^2z_{j}\right)^2\leq \frac{1}{2} E_h(0),
		\end{align*}
		and therefore
		$	-E_h(0)\leq L_{h,2}(t)|_0^T\leq E_h(0).$
		Next,  $I_4$ is simplified as the following
		\begin{align*}
			I_4=&-\frac{h}{8}\sum_{j=1}^{N}\int_{0}^{T}\left(\dot{z}_{j+\frac{1}{2}}+\dot{z}_{j-\frac{1}{2}}\right)\dot{z}_{j}dt-\frac{h}{8}\int_{0}^{T}\dot{z}_{N+\frac{1}{2}}\dot{z}_{N+1}dt\\
			=&-\frac{h}{8}\sum_{j=1}^{N}\int_{0}^{T}\dot{z}_{j+\frac{1}{2}}\dot{z}_{j}dt-\frac{h}{8}\sum_{j=0}^{N-1}\int_{0}^{T}\dot{z}_{j+\frac{1}{2}}\dot{z}_{j+1}dt-\frac{h}{8}\int_{0}^{T}\dot{z}_{N+\frac{1}{2}}\dot{z}_{N+1}dt\\
			=&-\frac{h}{8}\sum_{j=0}^{N}\int_{0}^{T}\dot{z}_{j+\frac{1}{2}}\dot{z}_{j}dt-\frac{h}{8}\sum_{j=0}^{N}\int_{0}^{T}\dot{z}_{j+\frac{1}{2}}\dot{z}_{j+1}dt\\
			=&-\frac{h}{4}\sum_{j=0}^{N}\int_{0}^{T}\dot{z}_{j+\frac{1}{2}}^2dt.
		\end{align*}
		Now, consider $I_5$ with \cref{poincare} and \cref{kdef} so that
		\begin{align*}
			I_5=&\frac{h}{4}\sum_{j=1}^{N}\int_{0}^{T}\delta_x^4z_jz_{j}dt-\frac{1}{4}\int_{0}^{T}\delta_x^3z_{N+\frac{1}{2}}z_{N+1}dt\\
			=&\frac{1}{4}\sum_{j=1}^{N}\int_{0}^{T}\left(\delta_x^3z_{j+\frac{1}{2}}-\delta_x^3z_{j-\frac{1}{2}}\right)z_{j}dt-\frac{1}{4}\int_{0}^{T}\delta_x^3z_{N+\frac{1}{2}}z_{N+1}dt\\
			=&\frac{1}{4}\sum_{j=1}^{N}\int_{0}^{T}\delta_x^3z_{j+\frac{1}{2}}z_{j}dt-\frac{1}{4}\sum_{j=0}^{N-1}\int_{0}^{T}\delta_x^3z_{j+\frac{1}{2}}z_{j+1}dt-\frac{1}{4}\int_{0}^{T}\delta_x^3z_{N+\frac{1}{2}}z_{N+1}dt\\
			=&\frac{1}{4}\sum_{j=0}^{N}\int_{0}^{T}\delta_x^3z_{j+\frac{1}{2}}z_{j}-\frac{1}{4}\sum_{j=0}^{N}\int_{0}^{T}\delta_x^3z_{j+\frac{1}{2}}z_{j+1}dt\\
			=&-\frac{h}{4}\sum_{j=0}^{N}\int_{0}^{T}\delta_x^3z_{j+\frac{1}{2}}\delta_xz_{j+\frac{1}{2}}dt\\
			=&-\frac{1}{4}\sum_{j=0}^{N}\int_{0}^{T}\left(\delta_x^2z_{j+1}-\delta_x^2z_{j}\right)\delta_xz_{j+\frac{1}{2}}dt\\
			=&-\frac{1}{4}\sum_{j=1}^{N+1}\int_{0}^{T}\delta_x^2z_{j}\delta_xz_{j-\frac{1}{2}}dt+\frac{1}{4}\sum_{j=0}^{N}\int_{0}^{T}\delta_x^2z_{j}\delta_xz_{j+\frac{1}{2}}dt\\
			=&\frac{h}{4}\sum_{j=0}^{N}\int_{0}^{T}\left(\delta_x^2z_{j}\right)^2dt-\underbrace{\frac{1}{4}\int_{0}^{T}\delta_x^2z_{N+1}\delta_xz_{N+\frac{1}{2}}dt}_{=0}\\
			\geq&\frac{h}{4}\sum_{j=0}^{N}\int_{0}^{T}\left(\delta_xz_{j+\frac{1}{2}}\right)^2dt\\
			=&\frac{h}{8}\sum_{j=0}^{N}\int_{0}^{T}\left(\delta_xz_{j+\frac{1}{2}}\right)^2dt+\frac{h}{8}\sum_{j=1}^{N+1}\int_{0}^{T}\left(\delta_xz_{j-\frac{1}{2}}\right)^2dt\\
			\geq&\frac{h}{16}\sum_{j=0}^{N}\int_{0}^{T}\left(\delta_xz_{j+\frac{1}{2}}+\delta_xz_{j-\frac{1}{2}}\right)^2dt\\
			=&\frac{1}{16h}\sum_{j=0}^{N}\int_{0}^{T}\left(y_j\right)^2dt\\
			=&\frac{1}{16h}\sum_{j=0}^{N}\int_{0}^{T}\left(\left(\frac{C}{P}-\frac{h^2}{4}\right)\left(\bm{A}_h\vec{k}\right)_j+k_j\right)^2dt\\
			=&\frac{1}{16h}\left(\frac{C}{P}-\frac{h^2}{4}\right)^2\sum_{j=0}^{N}\int_{0}^{T}\left(\left(\bm{A}_h\vec{k}\right)_j\right)^2dt+\frac{1}{8h}\left(\frac{C}{P}-\frac{h^2}{4}\right)\sum_{j=0}^{N}\int_{0}^{T}\left(\bm{A}_h\vec{k}\right)_jk_jdt+\frac{1}{16h}\sum_{j=0}^{N}\int_{0}^{T}k_j^2dt.
		\end{align*}
		Lastly, consider $I_6$ with  $z_N=\frac{z_{N+1}+z_{N-1}}{2}$
		\begin{align*}
			I_6&=-\frac{B}{8}\sum_{j=1}^{N}\int_{0}^{T}(\phi_{j+1}-\phi_{j-1})z_{j}dt+\frac{B}{4}\int_{0}^{T}\phi_Nz_{N+1}dt\\
			&=-\frac{B}{8}\sum_{j=1}^{N}\int_{0}^{T}\phi_{j+1}z_{j}dt+\frac{B}{8}\sum_{j=1}^{N}\int_{0}^{T}\phi_{j-1}z_{j}dt+\frac{B}{4}\int_{0}^{T}\phi_Nz_{N+1}dt\\
			&=-\frac{B}{8}\sum_{j=2}^{N+1}\int_{0}^{T}\phi_{j}z_{j-1}dt+\frac{B}{8}\sum_{j=0}^{N-1}\int_{0}^{T}\phi_{j}z_{j+1}dt+\frac{B}{4}\int_{0}^{T}\phi_Nz_{N+1}dt\\
			&=\frac{B}{8}\sum_{j=1}^{N}\int_{0}^{T}\phi_{j}y_jdt+\frac{B}{8}\int_{0}^{T}\phi_Nz_{N+1}dt-\frac{B}{8}\int_{0}^{T}\phi_Nz_{N}dt\\
			&=\frac{B}{8}\sum_{j=1}^{N}\int_{0}^{T}\phi_{j}y_jdt+\frac{B}{16}\int_{0}^{T}\phi_Ny_{N}dt.
		\end{align*}
		Therefore, the sum of  \eqref{firsttermssum}, $I_4$ and $L_{h,2}(t)|_0^T$ is estimated as the following
		\begin{equation}\label{sum}
			\begin{split}
				L_{h,1}(t)|_0^T+L_{h,2}(t)|_0^T+I_1+I_2+I_4\geq&-3E_h(0)-\frac{1}{2}\int_{0}^{T}\left(\dot{z}_{N+1}\right)^2dt\\
				&+\frac{h}{4}\sum_{j=0}^{N}\int_{0}^{T}\left(\dot{z}_{j+\frac{1}{2}}\right)^2dt+\frac{h}{2}\sum_{j=0}^{N}\int_{0}^{T}\left(\delta_x^2z_{j}\right)^2dt.
			\end{split}
		\end{equation}
		Moreover, the sum of $I_3$, $I_5$ and $I_6$ is estimated as
		\begin{equation}\label{shearsum}
			I_3+I_5+I_6\geq I_{7,1}+I_{7,2}+\frac{B}{8}\sum_{j=1}^{N}\int_{0}^{T}\phi_{j}y_jdt+\frac{5B}{16}\int_{0}^{T}\phi_Ny_{N}dt+\frac{1}{16h}\sum_{j=0}^{N}\int_{0}^{T}k_j^2dt.		
		\end{equation}
		where
		\begin{align*}
			I_{7,1}:=&\frac{1}{16h}\left(\frac{C}{P}-\frac{h^2}{4}\right)^2\sum_{j=0}^{N}\int_{0}^{T}\left(\left(\bm{A}_h\vec{k}\right)_j\right)^2dt+\frac{B^2}{8Ph}\left(\frac{C}{P}-\frac{h^2}{4}\right)\sum_{j=1}^{N}\int_{0}^{T}\left(\bm{A}_h\vec{k}\right)_{j}\left(\bm{A}_h\vec{k}\right)_{j+1}dt\\
			=&\frac{1}{32h}\left(\frac{C}{P}-\frac{h^2}{4}\right)^2\sum_{j=1}^{N}\int_{0}^{T}\left(\left(\bm{A}_h\vec{k}\right)_j\right)^2dt+\frac{1}{32h}\left(\frac{C}{P}-\frac{h^2}{4}\right)^2\sum_{j=1}^{N}\int_{0}^{T}\left(\left(\bm{A}_h\vec{k}\right)_{j-1}\right)^2dt\\
			&+\frac{B^2}{8Ph}\left(\frac{C}{P}-\frac{h^2}{4}\right)\sum_{j=1}^{N}\int_{0}^{T}\left(\bm{A}_h\vec{k}\right)_{j}\left(\bm{A}_h\vec{k}\right)_{j-1}dt\\
			&+\left(\frac{B^2}{8Ph}\left(\frac{C}{P}-\frac{h^2}{4}\right)+\frac{1}{32h}\left(\frac{C}{P}-\frac{h^2}{4}\right)^2\right)\int_{0}^{T}\left(\left(\bm{A}_h\vec{k}\right)_{N}\right)^2dt\\
			\geq&\left( \frac{1}{32h}\left(\frac{C}{P}-\frac{h^2}{4}\right)^2-  \frac{B^2}{16Ph}\left(\frac{C}{P}-\frac{h^2}{4}\right)\right)\sum_{j=1}^{N}\int_{0}^{T}\left(\left(\bm{A}_h\vec{k}\right)_j\right)^2dt\\
			&+\left( \frac{1}{32h}\left(\frac{C}{P}-\frac{h^2}{4}\right)^2 -\frac{B^2}{16Ph}\left(\frac{C}{P}-\frac{h^2}{4}\right) \right)\sum_{j=1}^{N}\int_{0}^{T}\left(\left(\bm{A}_h\vec{k}\right)_{j-1}\right)^2dt,
		\end{align*}
				and
		\begin{align*}
			I_{7,2}:=&\frac{1}{8h}\left(\frac{C}{P}-\frac{h^2}{4}\right)\sum_{j=1}^{N}\int_{0}^{T}\left(\bm{A}_h\vec{k}\right)_jk_jdt	+\frac{B^2}{8Ph}\sum_{j=1}^{N}\int_{0}^{T}\delta_xk_{j-\frac{1}{2}}\delta_xk_{j+\frac{1}{2}}dt\\
			=&\frac{1}{8h^2}\left(\frac{C}{P}-\frac{h^2}{4}\right)\sum_{j=1}^{N}\int_{0}^{T}\left(-\delta_xk_{j+\frac{1}{2}}+\delta_xk_{j-\frac{1}{2}}\right)k_jdt+\frac{B^2}{8Ph}\sum_{j=1}^{N}\int_{0}^{T}\delta_xk_{j-\frac{1}{2}}\delta_xk_{j+\frac{1}{2}}dt\\
			=&-\frac{1}{8h^2}\left(\frac{C}{P}-\frac{h^2}{4}\right)\sum_{j=1}^{N}\int_{0}^{T}\delta_xk_{j+\frac{1}{2}}k_jdt+\frac{1}{16h^2}\left(\frac{C}{P}-\frac{h^2}{4}\right)\sum_{j=0}^{N-1}\int_{0}^{T}\delta_xk_{j+\frac{1}{2}}k_{j+1}dt\\
			&+\frac{B^2}{8Ph}\sum_{j=1}^{N}\int_{0}^{T}\delta_xk_{j-\frac{1}{2}}\delta_xk_{j+\frac{1}{2}}dt\\
			=&-\frac{1}{8h^2}\left(\frac{C}{P}-\frac{h^2}{4}\right)\sum_{j=1}^{N}\int_{0}^{T}\delta_xk_{j+\frac{1}{2}}k_jdt+\frac{1}{16h^2}\left(\frac{C}{P}-\frac{h^2}{4}\right)\sum_{j=1}^{N}\int_{0}^{T}\delta_xk_{j+\frac{1}{2}}k_{j+1}dt\\
			&+\frac{B^2}{8Ph}\sum_{j=1}^{N}\int_{0}^{T}\delta_xk_{j-\frac{1}{2}}\delta_xk_{j+\frac{1}{2}}dt+\frac{1}{16h^2}\left(\frac{C}{P}-\frac{h^2}{4}\right)\int_{0}^{T}\delta_xk_{\frac{1}{2}}k_{1}dt\\
			=&\frac{1}{8h}\left(\frac{C}{P}-\frac{h^2}{4}\right)\sum_{j=1}^{N}\int_{0}^{T}\left(\delta_xk_{j+\frac{1}{2}}\right)^2dt\\
			&+\frac{B^2}{8Ph}\sum_{j=1}^{N}\int_{0}^{T}\delta_xk_{j-\frac{1}{2}}\delta_xk_{j+\frac{1}{2}}dt+\frac{1}{16h^3}\left(\frac{C}{P}-\frac{h^2}{4}\right)\int_{0}^{T}\left(k_{1}\right)^2dt\\
			=&\frac{1}{16h}\left(\frac{C}{P}-\frac{h^2}{4}\right)\sum_{j=1}^{N}\int_{0}^{T}\left(\delta_xk_{j+\frac{1}{2}}\right)^2dt+\frac{1}{16h}\left(\frac{C}{P}-\frac{h^2}{4}\right)\sum_{j=1}^{N}\int_{0}^{T}\left(\delta_xk_{j+\frac{1}{2}}\right)^2dt\\
			&+\frac{B^2}{8Ph}\sum_{j=1}^{N}\int_{0}^{T}\delta_xk_{j-\frac{1}{2}}\delta_xk_{j+\frac{1}{2}}dt+\frac{1}{16h}\left(\frac{C}{P}-\frac{h^2}{4}\right)\int_{0}^{T}\left(\delta_xk_{\frac{1}{2}}\right)^2dt\\
			=&\frac{1}{16h}\left(\frac{C}{P}-\frac{h^2}{4}\right)\sum_{j=1}^{N}\int_{0}^{T}\left(\delta_xk_{j+\frac{1}{2}}\right)^2dt+\frac{1}{16h}\left(\frac{C}{P}-\frac{h^2}{4}\right)\sum_{j=2}^{N+1}\int_{0}^{T}\left(\delta_xk_{j-\frac{1}{2}}\right)^2dt\\
			&+\frac{B^2}{8Ph}\sum_{j=1}^{N}\int_{0}^{T}\delta_xk_{j-\frac{1}{2}}\delta_xk_{j+\frac{1}{2}}dt+\frac{1}{16h}\left(\frac{C}{P}-\frac{h^2}{4}\right)\int_{0}^{T}\left(\delta_xk_{\frac{1}{2}}\right)^2dt\\
			=&\frac{1}{16h}\left(\frac{C}{P}-\frac{h^2}{4}\right)\sum_{j=1}^{N}\int_{0}^{T}\left(\delta_xk_{j+\frac{1}{2}}\right)^2dt+\frac{1}{16h}\left(\frac{C}{P}-\frac{h^2}{4}\right)\sum_{j=1}^{N}\int_{0}^{T}\left(\delta_xk_{j-\frac{1}{2}}\right)^2dt\\
			&+\frac{B^2}{8Ph}\sum_{j=1}^{N}\int_{0}^{T}\delta_xk_{j-\frac{1}{2}}\delta_xk_{j+\frac{1}{2}}dt+\frac{1}{16h}\left(\frac{C}{P}-\frac{h^2}{4}\right)\int_{0}^{T}\left(\delta_xk_{\frac{1}{2}}\right)^2dt\\
			\geq&\left( \frac{1}{16h}\left(\frac{C}{P}-\frac{h^2}{4}\right)-\frac{B^2}{16Ph} \right)\sum_{j=1}^{N}\int_{0}^{T}\left(\delta_xk_{j+\frac{1}{2}}\right)^2dt\\
			&+\left( \frac{1}{16h}\left(\frac{C}{P}-\frac{h^2}{4}\right)-\frac{B^2}{16Ph} \right)\sum_{j=1}^{N}\int_{0}^{T}\left(\delta_xk_{j-\frac{1}{2}}\right)^2dt.
		\end{align*}
		with the assumption
		\begin{equation}
		\label{rmk2}
		\frac{1}{2}\left(\frac{C}{P}-\frac{h^2}{4}\right)-\frac{B^2}{P} >0.
		\end{equation}
	 By \eqref{rmk2},  $I_{7,1}\geq0$ and $I_{7,2}\geq0.$ Thus,
		\begin{equation}\label{shearsum2}
			I_3+I_5+I_6\geq\frac{B}{8}\sum_{j=1}^{N}\int_{0}^{T}\phi_{j}y_jdt+\frac{5B}{16}\int_{0}^{T}\phi_Ny_{N}dt+\frac{1}{16h}\sum_{j=0}^{N}\int_{0}^{T}k_j^2dt.		
		\end{equation}
		Here, the positivity of the last two terms in \eqref{shearsum2} is significant for the rest of the proof. For this purpose,  consider two cases for the sign of $\frac{5B}{16}\int_{0}^{T}\phi_Ny_{N}dt$ in \eqref{shearsum2} as the following.
	\begin{itemize}
	\item [Case I:]	If \begin{equation}\label{onem1}
			\frac{5B}{16}\int_{0}^{T}\phi_Ny_{N}dt\geq0,
		\end{equation}
it is straightforward that
			\begin{equation}\label{positive}
			\frac{5B}{16}\int_{0}^{T}\phi_Ny_{N}dt+\frac{1}{16h}\int_{0}^{T}k_N^2dt\geq0.
		\end{equation}
	\item[Case II:]	Assume
	\begin{equation}
		\label{rmk1}
		\frac{2}{5}\left(\frac{C}{P}-\frac{h^2}{4}\right)-\frac{B^2}{P} >0.
		\end{equation}
If
		\begin{equation}\label{rmk5}
			\frac{5B}{16}\int_{0}^{T}\phi_Ny_{N}dt\leq0,
		\end{equation}
		\begin{align*}
			&\frac{5B}{16}\int_{0}^{T}\phi_Ny_{N}dt=	\frac{5B}{16}\int_{0}^{T}\frac{B}{2Ph}\left(\bm{A}_h\vec{k}\right)_N\left( \left(\frac{C}{P}-\frac{h^2}{4}\right) \left(\bm{A}_h\vec{k}\right)_N+k_N\right)dt\\
			&=\int_{0}^{T}\frac{5B^2}{32Ph}\left(\frac{C}{P}-\frac{h^2}{4}\right)\left(\bm{A}_h\vec{k}\right)_N^2+\frac{5B^2}{32Ph}\left(\bm{A}_h\vec{k}\right)_Nk_Ndt\\
			&=\frac{5B^2}{2P}\frac{1}{16h}\left[\int_{0}^{T}\left(\frac{C}{P}-\frac{h^2}{4}\right)\left(\bm{A}_h\vec{k}\right)_N^2+\left(\bm{A}_h\vec{k}\right)_Nk_Ndt\right]
		\end{align*}
Observe that the assumption \eqref{rmk5} directly implies
	$$\int_{0}^{T}\left(\frac{C}{P}-\frac{h^2}{4}\right)\left(\bm{A}_h\vec{k}\right)_N^2+\left(\bm{A}_h\vec{k}\right)_Nk_Ndt\leq 0, \quad \int_{0}^{T}\left(\bm{A}_h\vec{k}\right)_Nk_Ndt\leq 0.$$
Hence,
	\begin{align*}
					&\frac{5B}{16}\int_{0}^{T}\phi_Ny_{N}dt=\frac{5B^2}{2P}\frac{1}{16h}\left[\int_{0}^{T}\left(\frac{C}{P}-\frac{h^2}{4}\right)\left(\bm{A}_h\vec{k}\right)_N^2+\left(\bm{A}_h\vec{k}\right)_Nk_Ndt\right]\\
			&\underbrace{\geq}_{by ~\eqref{rmk}}\left(\frac{C}{P}-\frac{h^2}{4}\right)\frac{1}{16h}\left[\int_{0}^{T}\left(\frac{C}{P}-\frac{h^2}{4}\right)\left(\bm{A}_h\vec{k}\right)_N^2+\left(\bm{A}_h\vec{k}\right)_Nk_Ndt\right]\\
			&=\int_{0}^{T}\frac{1}{16h}\left(\frac{C}{P}-\frac{h^2}{4}\right)^2\left(\bm{A}_h\vec{k}\right)_N^2+\frac{1}{16h}\left(\frac{C}{P}-\frac{h^2}{4}\right)\left(\bm{A}_h\vec{k}\right)_Nk_Ndt\\
			&=\int_{0}^{T}\frac{1}{16h}\left(\frac{C}{P}-\frac{h^2}{4}\right)^2\left(\bm{A}_h\vec{k}\right)_N^2+\frac{2}{16h}\left(\frac{C}{P}-\frac{h^2}{4}\right)\left(\bm{A}_h\vec{k}\right)_Nk_Ndt\\
			&\qquad\underbrace{-\int_{0}^{T}\frac{1}{16h}\left(\frac{C}{P}-\frac{h^2}{4}\right)\left(\bm{A}_h\vec{k}\right)_Nk_Ndt}_{\geq0}\\
			&\geq\int_{0}^{T}\frac{1}{16h}\left(\frac{C}{P}-\frac{h^2}{4}\right)^2\left(\bm{A}_h\vec{k}\right)_N^2+\frac{2}{16h}\left(\frac{C}{P}-\frac{h^2}{4}\right)\left(\bm{A}_h\vec{k}\right)_Nk_Ndt\\
			&=\int_{0}^{T}\frac{1}{16h}\left(\frac{C}{P}-\frac{h^2}{4}\right)^2\left(\bm{A}_h\vec{k}\right)_N^2+\frac{2}{16h}\left(\frac{C}{P}-\frac{h^2}{4}\right)\left(\bm{A}_h\vec{k}\right)_Nk_Ndt\\
			&\quad+\frac{1}{16h}\int_{0}^{T}k_N^2dt-\frac{1}{16h}\int_{0}^{T}k_N^2dt\\
			&=\frac{1}{16h}\int_{0}^{T}\left( \left(\frac{C}{P}-\frac{h^2}{4}\right)\left(\bm{A}_h\vec{k}\right)_N+k_N \right)^2-\frac{1}{16h}\int_{0}^{T}k_N^2dt\\
			&=\frac{1}{16h}\int_{0}^{T}(y_N)^2dt-\frac{1}{16h}\int_{0}^{T}k_N^2dt,
		\end{align*}
which also implies \eqref{positive}.
		\end{itemize}
It is crucial to observe that the assumptions \eqref{rmk2} and \eqref{rmk1} lead to the assumption \eqref{rmk} in the statement of the theorem.
		Now, \eqref{positive} together with \eqref{onem1} and \eqref{positive} lead to
		\begin{align}
	\label{secondsum}				I_3+I_5+I_6\geq&\frac{B}{8}\sum_{j=1}^{N}\int_{0}^{T}\phi_{j}y_jdt+\frac{5B}{16}\int_{0}^{T}\phi_Ny_{N}dt+\frac{1}{16h}\sum_{j=0}^{N}\int_{0}^{T}k_j^2dt
			\geq\frac{B}{8}\sum_{j=1}^{N}\int_{0}^{T}\phi_{j}y_jdt.
		\end{align}
		Finally, by  the sum of \eqref{sum} and \eqref{secondsum},
			
		\begin{align*}
			0=&	\left[L_{h,1}(t)+L_{h,2}(t)\right]|_0^T+I_1+I_2+I_3+I_4+I_5+I_6\\
			\geq&-3E_h(0)-\frac{1}{2}\int_{0}^{T}\left(\dot{z}_{N+1}\right)^2dt+\underbrace{\frac{h}{4}\sum_{j=0}^{N}\int_{0}^{T}\left(\dot{z}_{j+\frac{1}{2}}\right)^2dt+\frac{h}{2}\sum_{j=0}^{N}\int_{0}^{T}\left(\delta_x^2z_{j}\right)^2dt+\frac{B}{8}\sum_{j=1}^{N}\int_{0}^{T}\phi_{j}y_jdt}_{\geq\frac{T}{2}E_h(t)}\\
			\geq&\left(\frac{T}{2}-3\right)E_h(0)-\frac{1}{2}\int_{0}^{T}\left(\dot{z}_{N+1}\right)^2dt.
		\end{align*}\normalsize
		This leads to the desired observability inequality  \eqref{MAIN} in the statement of the theorem.
	\end{proof}

\section{Simulations by Numerical Experiments}

\color{black}

First, consider a three-layer beam where the top outer layer is piezoelectric (Lead Zirconate Titanate - PZT), the bottom outer layer is aluminum, and the viscoelastic layer is silicone rubber, for which appropriate material constants are given in Table \ref{constants}.
	\begin{table}[htb!]		
	\centering
	\begin{tabular}{|c|c|c|c|}
		\hline
		Parameter& PZT& Rubber & Aluminum\\
		\hline
		Mass density (kg/$\text m^3$) & $\rho_{1}=7500  $ & $\rho_{2 }= 1250$ &$ \rho_{3}= 2710$ \\
		\hline
		Thickness (m) & $h_{1}=0.01$ & $h_{2 }=0.03$ &$ h_{3}=0.01$ \\
		\hline
		Young's Modulus (GPa) & $E_{1}=72$ & $E_{2}=0.05$ & $E_{3}=70$\\
		\hline
		Shear Modulus (GPa) &  $G_{1}=27$ & $G_{2}=0.01$ & $G_{3}=25$ \\
		\hline
		Poisson's Ratio & $\nu_{1}=0.31$ & $\nu_{2}=0.47$ &  $\nu_{3}=0.33$\\
		\hline
	\end{tabular}
	\caption{Material constants for each layer, i.e. see \url{https://www.matweb.com}.}
	\label{constants}
\end{table}

For  $D_i=\frac{E_i}{12(1-\nu_i^2)}$ with $D_ih_i^3$ being the modulus of flexural rigidity of each layer $i=1,2,3$, the positive constants $B,C,P$ for the three-layer beam model \eqref{eq1a}  are  as the following
\begin{align*}
	B=\frac{G_2(h_1+2h_2+h_3)}{2h_2(D_1h_1^3+D_3h_3^3)},&& C=\frac{G_2}{h_2(D_1h_1^3+D_3h_3^3)},&&P=\frac{G_2(D_1h_1+D_3h_3)}{12h_2^2D_1D_3h_1h_3(D_1h_1^3+D_3h_3^3)},
\end{align*}
For the material constants of the materials given in Table \ref{constants}, $B\approx1011.318, C\approx 25282.944,  P\approx 2130860.555.$

Letting $N=20$ and defining the mesh size $h:=\frac{1}{21}$, consider a uniform discretization of the interval $[0,L]$: $0=x_0<x_1<...<x_i=i*h<\ldots <x_{20}<x_{21}=L=1.$
Let $\left\{(z_i,\phi_i)(t)\approx (z,\phi)(x_i,t)\right\}_{i=0}^{N+1}$  be the approximated solutions in \eqref{clampeddiscrete} of the solutions $(z,\phi)(x,t)$ of \eqref{eq1a}-\eqref{eq1d}.
To show the strength of the proposed model reduction \eqref{clampeddiscrete} and the corresponding sensor design, consider the standard Finite Difference approximations of  \eqref{raymulti}. Notice that a standard Finite Difference-based model reduction without any averaging or middle nodes can be considered as the following
\begin{equation}\label{classicfd}
	\begin{cases}
		\ddot{z}+\delta_x^4z_i-\frac{B}{h}\left(\phi_{i+1}-\phi_{i}\right)=0,&\\
		-C\left(\frac{\phi_{i+1}-2\phi_i+\phi_{i-1}}{h^2}\right)+P\phi_i+\frac{B}{2}\left(\frac{z_{i+2}-3z_{i+1}+3z_{i}-z_{i-1}}{h^3}\right)=0,& i=1,...,N,\\
		\frac{z_{N+2}-3z_N+3z_{N-1}-z_{N-2}}{h^3}-B\phi_{N}=0,\\
		\phi_0=0, ~\phi_{N+1}=\phi_N, ~ z_0=z_{-1}=0, ~z_{N+2}=2z_{N+1}-z_N,\\
		z_i(0)=z^0(x_i),\quad \dot{z}_i(0)=z^1(x_i),& i=0,1,\ldots, N+1,
	\end{cases}
\end{equation}
where the boundary equation does not have its inertia term anymore. From now on, model reductions are shown with their abbreviations, i.e.  Finite Differences (FD) and Order-Reduced Finite Differences (ORFD).

First, the eigenvalues of the systems FD and ORFD are purely imaginary. The very immediate observation is that the gap between the two consecutive high-frequency eigenvalues, i.e. $|\lambda_{k+1}-\lambda_k|,$ for the FD model tend to zero which generally causes the lack of observability of the FD model reduction, see Fig. \ref{eigvalsconsfig}. A similar phenomenon is recently observed for an FD model  reduction yet with fully-hinged boundary conditions \cite{Aydin-O1}. As proved in Section 5,  the ORFD model reduction does not have any issue of this sort.

\begin{figure}[htb!]
	\centering
	\includegraphics[width=0.6\linewidth]{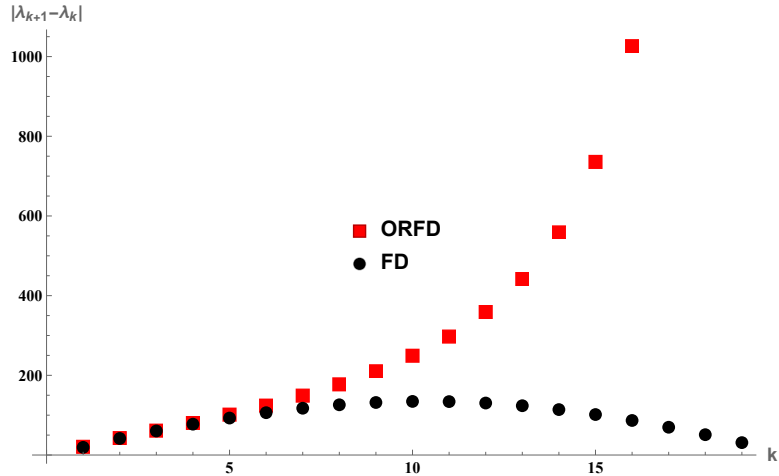}
	\caption{The uniform gap decreases for the (FD) model. }
	\label{eigvalsconsfig}
\end{figure}

To show the strength of the proposed sensor/actuator design for the closed-loop model, consider the PDE model \eqref{raymulti} with  the stabilizing feedback boundary condition  $z'''(L,t) - B\phi (L,t)=\xi \dot z(L,t)$ with $\xi>0.$ The closed-loop system is now dissipative and  is shown to be exponentially stable \cite{IFAC}.
Therefore, the homogeneous transverse shear boundary conditions in the ORFD model \eqref{clampeddiscrete} and the FD model \eqref{classicfd} are changed to the following
\begin{equation}\label{discretefd}
\left\{	\begin{array}{ll}
	\frac{z_{N+2}-3z_N+3z_{N-1}-z_{N-2}}{h^3}-\frac{h}{4}\left(\ddot{z}_{N+1}+\ddot{z}_{N}\right)-B\phi_{N}=\xi \dot{z}_{N+1},& {\rm (ORFD) ~model}\\
		\frac{z_{N+2}-3z_N+3z_{N-1}-z_{N-2}}{h^3}-B\phi_{N}=\xi \dot{z}_{N+1},& {\rm (FD) ~model}.
	\end{array}
\right.
\end{equation}

For the simulations below, the original model \eqref{eq1a} is re-scaled for the time variable, 
  i.e. $t^*=\sqrt{\frac{m}{A}}t$ and $\frac{d}{dt}\to \sqrt{\frac{A}{m}}\frac{d}{dt*}$. For the material constants given in Table \ref{constants}, $\sqrt{\frac{m}{A}}\approx0.1$. The initial conditions are chosen to be "box-like" as $z_i(0)=\dot{z}_i(0)=10^{-3}\cdot \chi_{[0.25,0.75]}(x_i)$ for $i=1,2,...,N+1$.

For the feedback gain $\xi=5,$ it is observed in Fig. \ref{bendingfig} that the bending solution $z(x,t)$ of the FD model reduction \eqref{classicfd},\eqref{discretefd}, the high-frequency vibrations still remain as the time evolves for $T=10$ sec., which indicates a lack of exponential stability uniformly.  However, the bending solution $z(x,t)$ of the ORFD model reduction \eqref{clampeddiscrete},\eqref{discretefd} shows rapid stabilization.
\begin{figure}[htb!]
	\centering
	\includegraphics[width=\linewidth]{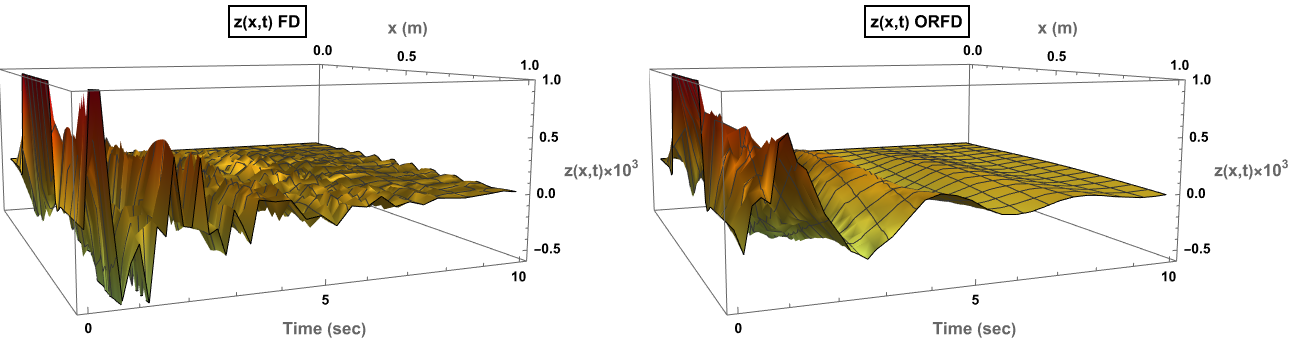}
	\caption{Overall bending of the beam $z(x,t)$ for FD vs. ORFD.}
	\label{bendingfig}
\end{figure}

As seen in Figure \ref{eigvalsfig}, the eigenvalue distributions of both models show a substantial difference. The ones for the model reduction \eqref{classicfd},\eqref{discretefd} have the real part of the eigenvalues  approach the imaginary axis as $h\to 0,$ which indicates the lack of exponential stability,  whereas the ones of the ORFD model reduction \eqref{clampeddiscrete},\eqref{discretefd} are strictly bounded away from the imaginary axis. The normalized energies of solutions and feedback sensor data for each model reduction confirm the same behavior, see Figures \ref{sensorfig} and \ref{energyfig}.
\begin{figure}[htb!]
	\centering
	\includegraphics[width=0.6\linewidth]{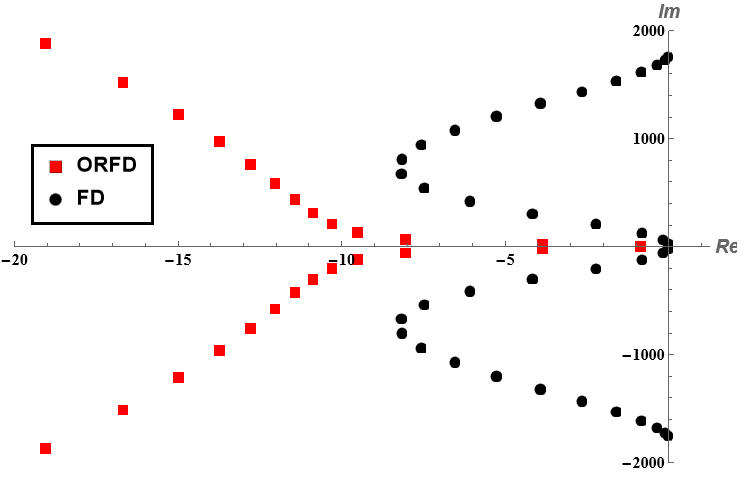}
	\caption{The real parts of the eigenvalues of the FD model reduction converge to the imaginary axis. This is not the case for the (ORFD) model reduction.}
	\label{eigvalsfig}
\end{figure}
\begin{figure}[htb!]
	\centering
	\includegraphics[width=0.8\linewidth]{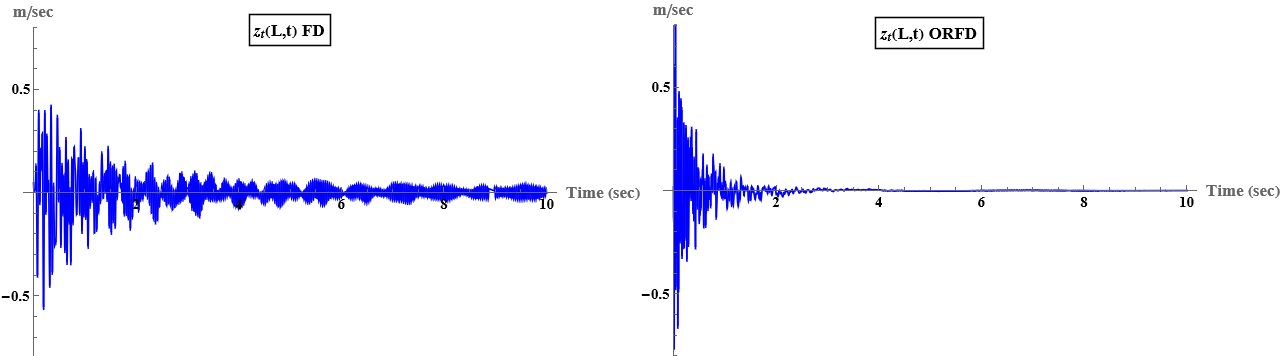}
	\caption{Boundary sensor measurements $\dot{z}(L,t)$ for FD and ORFD.}
	\label{sensorfig}
\end{figure}
\begin{figure}[htb!]
	\centering
	\includegraphics[width=0.6\linewidth]{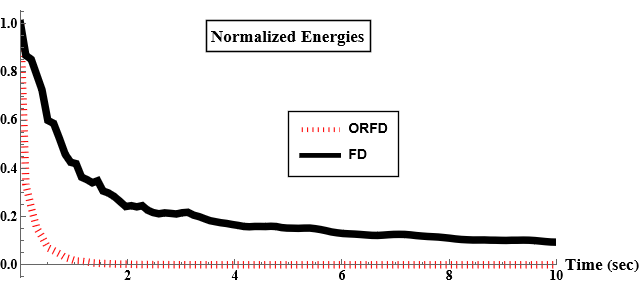}
	\caption{Normalized energies the FD and ORFD model reductions. The FD model shows a lack of exponential stability.}
	\label{energyfig}
\end{figure}

\color{black}
	\section{Conclusions, Ongoing Work, and Future Directions}

First, it is shown that  the PDE model \eqref{eq1a} - \eqref{eq1d}
is exactly observable with a single boundary sensor design by the so-called ``multipliers'' method.  A first-order model reduction is considered by a novel Finite-Difference approximation on equidistant grid points and average operators. This type of model reduction allows the approximations of the bending vibrations and shear angles of the beam for ``in-between nodes''. Therefore, a new variation of the boundary equation is obtained as a second-order ODE. The discretized energy  is defined and shown to be conservative for all time. Moreover, the approximated model is shown to retain uniform observability, as the mesh parameter $h\to 0$, without the need of spurious high-frequency filtering. The main hurdle here is the absorption of the coupling of the shear with bending dynamics sharply different from a single-layer Euler-Bernoulli beam \cite{Guo1,Guo2,L-Z}.

The immediate extension of the current work, also currently under investigation \cite{Aydin-O3},  is to prove the  controllability/exponential stability of reduced model of \eqref{raymulti} uniformly as $h\to 0$  with a single boundary controller. The proof requires a more delicate analysis due to  the matrix coefficients $\bf{B,C,P}$ in the PDE  model \eqref{raymulti} in comparison to the $B,C,P$ in \eqref{eq1a}. The corresponding Wolfram Demonstrations Projects are currently at the stage of submission \cite{WDP4,WDP5}. 	

 One other obvious extension of these results, which is also under investigation, is to improve the sub-optimal observation times $T>2$ and $T>6$ in the PDE model (Theorem \ref{cf-imp1}) and its model reduction (Theorem \ref{mainthm}), repsectively. In fact, exploring other multipliers for getting rid of the large shear condition \eqref{rmk} in Theorem \ref{mainthm}, arising due to the use of discrete multipliers method, is strongly needed since this condition does not exist for the PDE model (Theorem \ref{cf-imp1}).

	A  natural future direction is to prove  the uniform controllability and uniform exponential stability result (as $h\to 0$) for a reduced-order model of the Rao-Nakra-type sandwich beam model \cite{O3,O-IEEE2}, where the longitudinal dynamics for the outer layers are retained in the model. The analysis here together with the approach in \cite{ozer_uniform_2019} are  directly applicable to get the desired results.

\end{document}